\documentclass[12pt,a4paper]{article}

\usepackage[english]{babel}
\usepackage{epsfig}
\usepackage{amsmath}
\usepackage{amssymb}
\usepackage{amsbsy}
\usepackage{amsfonts}
\usepackage{mathrsfs} 
\usepackage{array}
\usepackage{verbatim}
\usepackage{graphicx}
\usepackage{color}
\usepackage{tabularx}
\usepackage{bm}
\usepackage{booktabs}
\usepackage[font=small]{caption}
\usepackage{subfig}
\usepackage{stmaryrd}
\usepackage{textcomp}

\begin{document}

\title{A seamless extension of DG methods for hyperbolic problems to unbounded domains}
\author{Tommaso Benacchio $^{(1)}$, Luca Bonaventura $^{(2)}$}
\maketitle

\begin{center}
 
{\small
$^{(1)}$  Met Office\\
 Exeter \\
United Kingdom\\
{\tt tommaso.benacchio@metoffice.gov.uk}
}
\vskip 0.5cm
{\small
$^{(2)}$ MOX -- Modelling and Scientific Computing, \\
Dipartimento di Matematica, Politecnico di Milano \\
Via Bonardi 9, 20133 Milano, Italy\\
{\tt luca.bonaventura@polimi.it}
}

\end{center}

\date{}


\vspace*{0.5cm}


\vspace*{0.5cm}

\pagebreak

\abstract{We consider spectral discretizations of hyperbolic problems on unbounded domains using Laguerre basis functions. Taking as  model problem  the scalar
advection equation, we perform a comprehensive stability analysis that includes strong collocation formulations, nodal and modal weak formulations, with either inflow or outflow boundary conditions, using either Gauss - Laguerre or Gauss - Laguerre - Radau quadrature nodes and based on either  scaled Laguerre functions or scaled Laguerre polynomials. We show that some of these combinations give rise to intrinsically unstable schemes, while
the combination of  scaled Laguerre functions with Gauss - Laguerre - Radau nodes appears to be stable for both
strong and weak formulations.
 We then show how  a modal discretization approach for hyperbolic systems on an unbounded domain can be naturally and seamlessly coupled to a discontinuous finite element discretization on a finite domain. 
 Examples of one dimensional hyperbolic systems are solved with the proposed domain decomposition technique.
  The errors obtained with the proposed approach are found to be small, enabling the 
use of the coupled scheme for the simulation of   Rayleigh damping layers in the semi-infinite part. Energy errors and reflection ratios of the scheme in absorbing 
wavetrains and single Gaussian signals show that a small number of nodes in the semi-infinite domain are sufficient to damp the waves. The theoretical insight and numerical results corroborate 
previous findings by the authors and establish the scaled Laguerre functions-based discretization as a flexible and efficient tool for absorbing layers  as well 
as for the accurate simulation of waves in unbounded regions.}

\pagebreak
 
\section{Introduction}
\label{intro} \indent

Accurate numerical solution of wave propagation problems in unbounded domains remains to the present day an unsolved challenge. Practical applications of great importance include
 modelling of the upper terrestrial atmosphere and of the solar corona, space weather simulations and propagation of electromagnetic waves from a localized source into the far field.
In other contexts such as numerical weather prediction, for computational reasons the domain of the hyperbolic differential problem is restricted to a bounded region of interest, the
size of which depends on the phenomena under consideration. Boundary conditions are then imposed to the new domain such
that outgoing waves can propagate from the bounded domain without spurious reflections. The issue of how to set these open boundary conditions in an analytically consistent and
numerically accurate and efficient way has been the object of research for the past four decades \cite{givoli:2008, rasch:1986,tsynkov:1998}. 

On the one hand, analytical approaches have aimed at imposing radiation or characteristic boundary conditions modelling outgoing
disturbances as solutions of reduced model dynamics, see \cite{engquist:1977, israeli:1981} as well as 
\cite{dea:2009, dea:2011, neta:2008} more specifically on environmental fluid dynamics.
On the other hand, tackling the problem numerically requires the placement of a buffer region next to the artificial boundary where the outgoing disturbances are relaxed towards a prescribed
external solution, commonly the absence of perturbations, that in turn calls for the addition of a diffusive or reactive term in the buffer region
\cite{durran:1983,  klemp:2008, klemp:1978, lavelle:2008}. Especially challenging in the setup of these absorbing layers is the choice of resolution, as finer grid spacings enable a better
absorption of the outgoing waves, yet incur a higher computational cost. Advanced perfectly matched layer formulations, popular in electro- and elastodynamics as well as computational
aeroacoustics, aim at matching analytical formulations on either side of the interface between the bounded and unbounded region, while optimally tuning relaxation parameters in the
layer \cite{abarbanel:1999, berenger:1994, hesthaven:1998, modave:2017, navon:2004,  rabinovich:2015, yang:1999}. Despite the substantial efforts in this context, the choice of parameters such
as layer thickness and relaxation coefficients is still largely driven by bespoke criteria \cite{modave:2010}. Yet another approach is given by infinite elements, whose shape functions
mimick the asymptotic behaviour of the solution at infinity \cite{astley:2000, gerdes:2000}. Use of infinite elements-based methodologies does not require splitting the original unbounded
domain, yet it entails mapping and modulation by decay functions.

In most cases, both the coupling conditions and the damping setup often rest on restrictive assumptions on the form of outgoing waves \cite{appelo:2009}.
However, approaches demanding an \emph{ad hoc} identification of waves appear ill-suited to the manifold nature of wavelike solutions impinging on the artificial ceiling of the computational
domain in environmental fluid dynamics problems \cite{durran:2015}. Furthermore, as more computational power becomes available,
modellers are faced with the complementary requirements of simulating ever larger portions of unbounded domains while achieving efficient control over spurious reflections at
the upper boundary \cite{giorgetta:2006, horinouchi:2003}. As a point of reference, numerical models for operational weather forecasting and climate research routinely devote up to
a quarter of the computational power in their simulations to eliminate unwanted perturbations in the upper atmosphere regions. This cost is bound to increase in future high-resolution
models as the size needed for the absorbing layer increases with increasing horizontal grid-spacing.

In order to overcome the drawbacks of currently employed approaches, in \cite{benacchio:2010, benacchio:2013} the authors devised the first application of scaled Laguerre function discretizations
to wave propagation problems, using a spectral collocation method for the shallow water equations on the positive half line. By tuning the scaling parameter, the set of quadrature nodes associated with the basis spans a differently-sized portion of the unbounded domain
\cite{ shen:2001, shen:2009,   wang:2009, zhuang:2010}. As shown in \cite{benacchio:2013}, this approach to unbounded domains could be coupled to standard discretizations on finite size domains, in order to
achieve a convenient framework for effective and economical implementation of open boundary conditions. The strategy entailed a significant reduction in
the computational cost of an absorbing layer.

In the present work, we first assess the many possible variants for a discretization of hyperbolic wave propagation
 problems on unbounded domains using Laguerre basis functions. Choosing the linear advection
 equation on the domain $ {\mathbb R}^+=[0,+\infty)$ as model problem,
 we find that discretizations   based on scaled Laguerre functions, rather than scaled Laguerre polynomials, and Gauss-Laguerre-Radau quadrature nodes, rather than Gauss-Laguerre nodes, yield stable spectra
 for the semidiscrete problem in space. 
%
%
%
%
%

We then extend the framework of \cite{benacchio:2013} and couple a Laguerre spectral   discretization  to a finite element discretization. More specifically, the
domain is split into a bounded interval of size $L$, where the system is discretized with a standard Discontinuous Galerkin (DG) approach, and a complementary unbounded interval, where the chosen
Laguerre spectral approach is employed, so that  ${\mathbb R}^+=[0,L]\cup[L,+\infty)$. 
We focus here on DG discretization approaches, since the discontinuous nature of the basis functions allows for easier coupling to the discretization
on the semi-infinite domain. We show that a modal Laguerre spectral element scheme can be coupled seamlessly to a modal DG discretization on the finite domain, providing a convenient
framework to extend finite domain discretizations and to implement open boundary conditions efficiently via absorbing layers. Numerical experiments are then performed, repeating some of the tests already presented
in  \cite{benacchio:2013}, showing that the accuracy of the proposed approach and its efficiency in implementing
absorbing layers improve those of our previous work.
%
%
%
%
 
The structure of the paper is as follows. The spectral stability analysis of Laguerre discretizations of scalar advection is discussed
in Section \ref{analysis}. The case of hyperbolic systems is tackled in Section \ref{modal_sys},
while in Section  \ref{coupled}  an outline of the coupling strategy is presented. Section \ref{tests} contains the numerical results and
Section \ref{conclu} some conclusions and discussion of future work.
  
\section{Stability analysis of scaled Laguerre discretizations of the advection equation}
\label{analysis}
We plan to analyze different discretizations of hyperbolic equations
on semi-infinite domains,
in order to understand which approach is most convenient to couple
with finite element discretizations on finite domains.
For this purpose, we start considering the advection equation as a prototype of hyperbolic problems.
In strong form, the advection equation is given by  
\begin{equation}
\frac{\partial q}{\partial t} + u\frac{\partial q}{\partial z} =0 \qquad z \in{\mathbb R}^+=[0,+\infty).  
\label{advection}
\end{equation}
This equation should be complemented by the condition
$$\lim_{z\rightarrow+\infty}q(z,t)=0$$
and by Dirichlet boundary condition $q(0,t)=q_L$ when $u>0$,
while  in the case $u<0$ no boundary conditions should be assigned,
nor required by the discretization method. It is also assumed that

To derive    the corresponding weak form, we integrate equation \eqref{advection}
against a test function $  \varphi $ to obtain
\begin{equation}
 \frac{d}{dt}\int_0^{+\infty} \varphi(z) q(z,t)\omega(z)\,dz
+u \int_0^{+\infty} \varphi(z) \frac{\partial q }{\partial z}(z,t)\omega(z)\,dz =0
\label{weak_form}
\end{equation}
If scaled Laguerre functions are used, $\omega=1$. We can then integrate by parts
and obtain, in the inflow case,
\begin{equation}
 \frac {d }{d t}\int_0^{+\infty} \varphi(z) q(z,t) \,dz - uq_L\varphi(0) 
-u \int_0^{+\infty} \varphi^{\prime}(z) q (z,t)\,dz =0.
\label{weak_form_ibp_in}
\end{equation}
In the outflow case, in order to have a well posed problem, the solution should not
be prescribed at the boundary, so that one should write instead
\begin{equation}
 \frac{d}{dt}\int_0^{+\infty} \varphi(z)  q(z,t)\,dz - uq(0,t)\varphi(0) 
-u \int_0^{+\infty} \varphi^{\prime}(z) q (z,t)\,dz =0.
\label{weak_form_ibp_out}
\end{equation}
In case we want to use scaled Laguerre polynomials as basis and test functions, we have to  
 assume instead  $\omega(z) =\exp{(-\beta z)}$ in  \eqref{weak_form}
before integrating by parts. Since
\begin{equation}
  \varphi\frac{\partial q}{\partial z}\omega = \frac{\partial }{\partial z}\left( \varphi q\omega \right) 
-\varphi q\frac{\partial \omega}{\partial z} -\frac{\partial  \varphi }{\partial z} q\omega 
= \frac{\partial }{\partial z}\left(\varphi q \omega \right) 
+\beta \varphi q\omega -\frac{\partial\varphi}{\partial z} q\omega,
\end{equation}
where we have used the fact that $\frac{\partial \omega}{\partial z}=-\beta \omega$,
we can obtain, in the inflow case,
\begin{multline}
\frac {d }{d t}\int_0^{+\infty} \varphi(z) q(z,t) \omega(z)\,dz +\beta u\int_0^{+\infty} \varphi(z) q(z,t) \omega(z)\,dz \\ 
- uq_L\varphi(0) \omega(0) -u \int_0^{+\infty} \varphi^{\prime}(z) q(z,t) \omega(z)\,dz =0 \label{weak_form_poly_ibp_in}
\end{multline}
and in the outflow case
\begin{multline}
\frac {d}{dt}\int_0^{+\infty} \varphi(z)  q(z,t) \omega(z)\,dz + \beta u\int_0^{+\infty} \varphi(z)  q(z,t) \omega(z)\,dz \\ 
- uq(0,t)\varphi(0) \omega(0) - u \int_0^{+\infty} \varphi^{\prime}(z) q (z,t) \omega(z)\,dz =0. \label{weak_form_poly_ibp_out}
\end{multline}

Several   alternatives are then possible for the numerical discretization.
As test and basis functions, one may choose
a) scaled Laguerre functions or b) scaled  Laguerre polynomials.
Option a) has the advantage of avoiding the presence of a weight in the scalar
product, while option b) has the advantage of allowing to approximate
constant functions on ${\mathbb R}^+$.
Equation \eqref{advection} can then be discretized   1) in strong form by a collocation approach
using Gauss-Laguerre Radau (GLR) nodes, which is the only practical alternative if Dirichlet b.c. have to be imposed,
2) in weak form, using either GLR or Gauss Laguerre (GL) nodes for numerical integration.
Furthermore for option 2), either a nodal  discretization (2n) can be used, in which
the basis functions  considered are Lagrange basis functions associated with the chosen
integration nodes, or a modal discretization (2m).
In all cases, spatial discretization results in an ODE system
\begin{equation}
\frac{d \mathbf{q}}{dt} = \mathbf{A} \mathbf{q}+ \mathbf{g},
\label{disc_advection}
\end{equation}
where $\mathbf{q}$ denotes the vector of the degrees of freedom of the spatially discretized problem,
$\mathbf{A}$ the discretization of the advection operator and  $\mathbf{g}$
is a source term that depends on the Dirichlet b.c. in the inflow boundary $u>0$ case,
while $\mathbf{g}=\mathbf{0}$ in the outflow boundary $u<0$ case.

We want to study the eigenvalue structure of $\mathbf{A}$ in order to better understand possible
stability and accuracy issues in the coupling to discretizations on finite domains. 
We adopt mostly the notation of \cite{benacchio:2010}, \cite{benacchio:2013} and  
denote by $\beta$ a scaling parameter,  by $z^{\beta}_l,\,l=0,\dots,M$ the scaled
GLR (SGLR) or the scaled GL (SGL) nodes. We denote by
$\omega^{\beta}_l$ the associated weights and
$\hat h^{\beta}_l$ are the associated Lagrange interpolation functions. 

\subsection*{Strong form, collocation approach with SGLR  nodes and scaled Laguerre functions}

We first consider the case a1   (strong collocation form with SGLR  nodes and scaled Laguerre
functions). We consider the outflow and inflow cases separately.
For the inflow case, 
\begin{equation}
q_i^{\prime} = -uq_L (\hat h^{\beta}_0)^{\prime}(z_i^{\beta}) - u\sum_{j=1}^M q_j (\hat h^{\beta}_j)^{\prime}(z_i^{\beta}) \qquad i=1,\dots,M.
\label{coll_sglr_in}
\end{equation}
In the outflow case, one has instead
\begin{equation}
q_i^{\prime} =   -u\sum_{j=0}^M q_j (\hat h^{\beta}_j)^{\prime}(z_i^{\beta})
  \qquad i=0,\dots,M.
\label{coll_sglr_out}
\end{equation}
\begin{figure}[htb]
\centering
\subfloat[]{\includegraphics[scale=0.45]{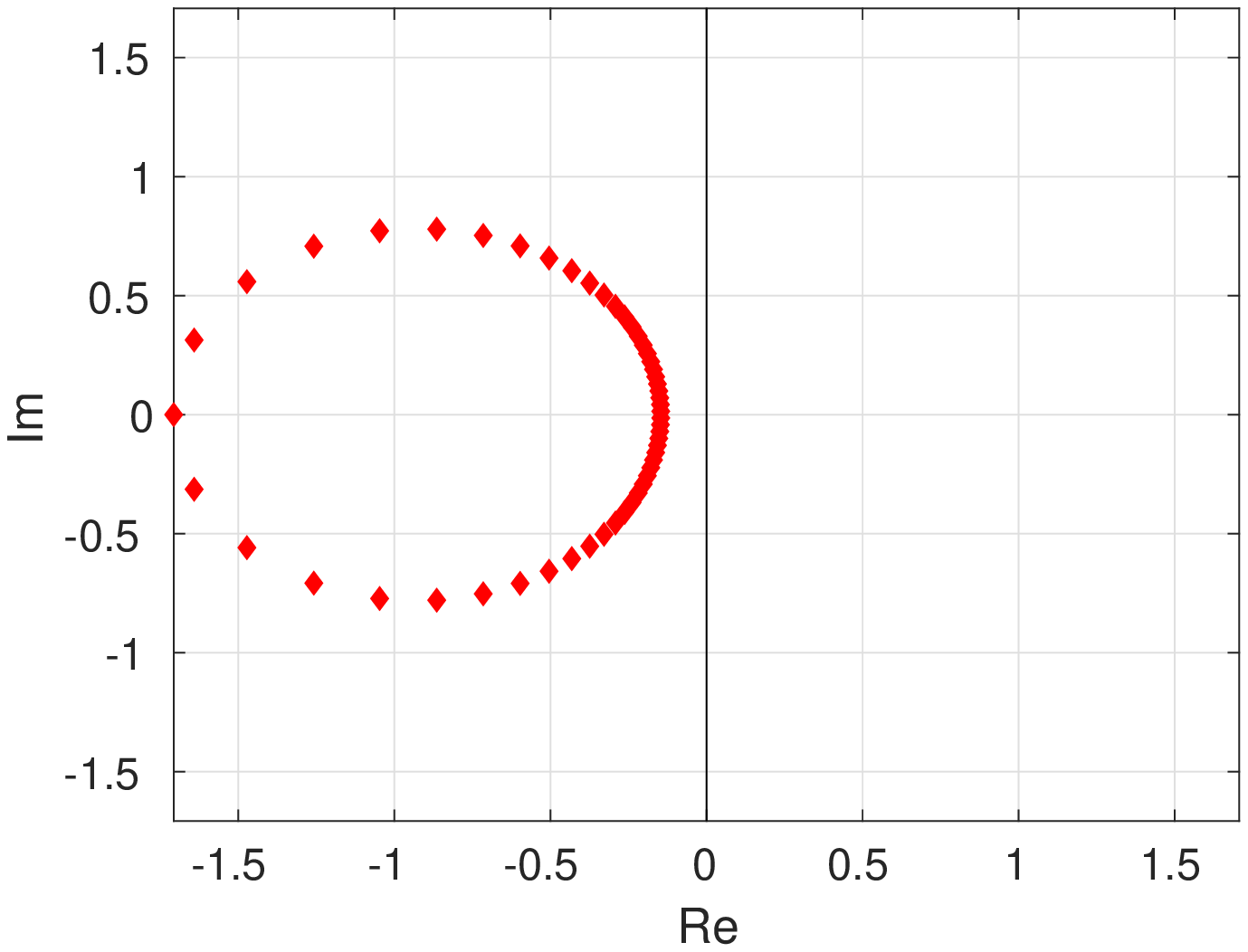}}\quad
\subfloat[]{\includegraphics[scale=0.45]{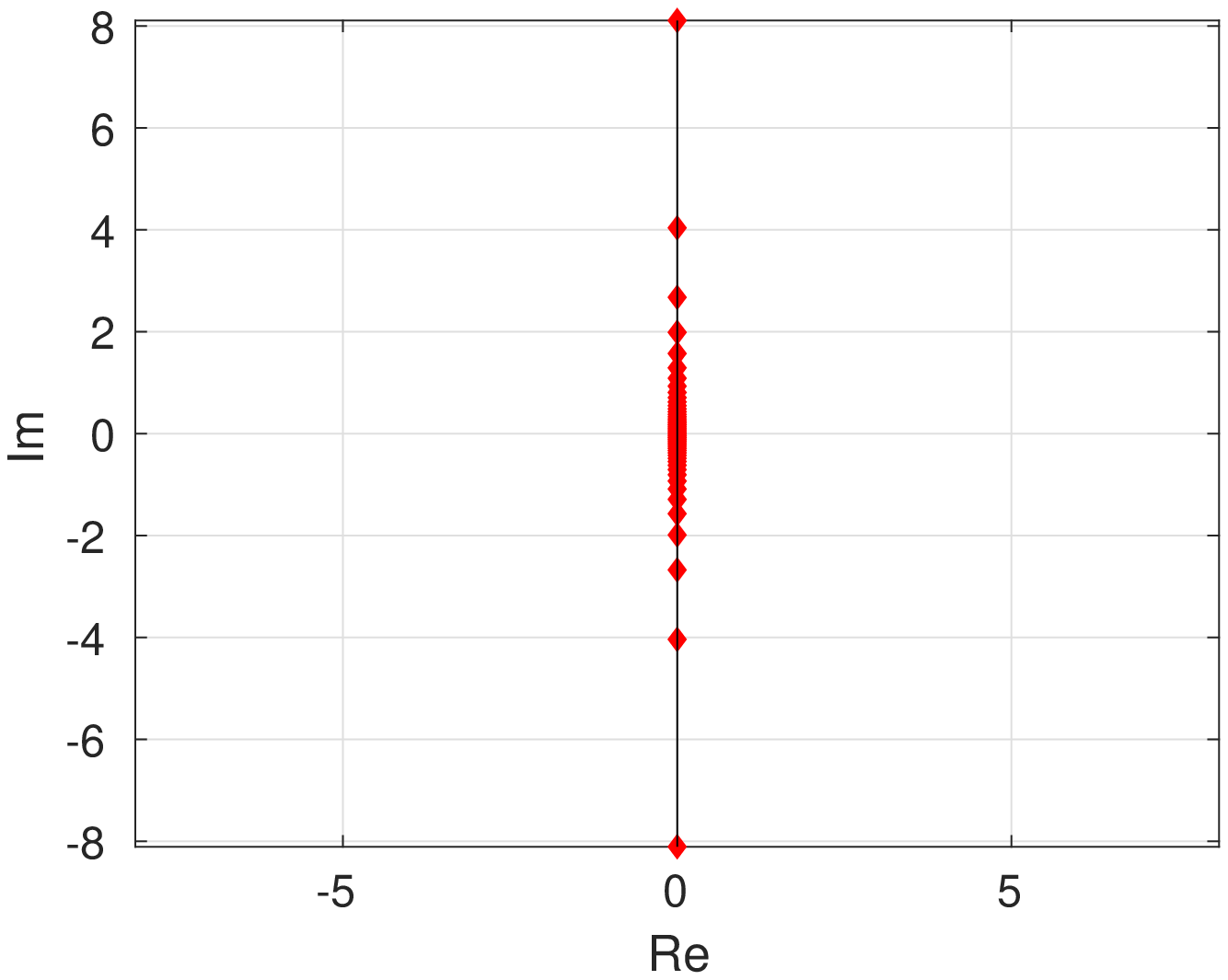}}
\caption{Eigenvalues in the a1 case, for parameter values  $\beta=1$, $M=50$
in the a) outflow  case $u=-1$ and b) inflow  case $u=1$. }
\label{fig:a1_m50_beta1}
\end{figure}
In matrix notation, equation \eqref{coll_sglr_in} yields in the inflow case 
an ODE system like equation \eqref{disc_advection} with
$\mathbf{q}=[q_1,\dots,q_M]^T$, 

\begin{equation} 
\mathbf{A}=-u(\mathbf{D}_{\beta})_M, \qquad \mathbf{g}=-uq_L\left[(\hat h^{\beta}_0)^{\prime}(z_1^{\beta}),
\dots,(\hat h^{\beta}_0)^{\prime}(z_M^{\beta})\right] ^T
\end{equation}
where, in the notation of \cite{benacchio:2013}, $\mathbf{D}_{\beta}$ denotes the SLGR differentiation matrix and
$(\mathbf{D}_{\beta})_M$ denotes the $M\times M$ matrix obtained selecting the last $M$ rows and columns of $\mathbf{D}_{\beta}$.
For the outflow case, $\mathbf{q}=[q_0,\dots,q_M]^T$, $\mathbf{g}=\mathbf{0}$ and $\mathbf{A}=-u\mathbf{D}_{\beta}$. Examples of eigenvalue
plots for $\mathbf{A}$ in the a1 case are shown in Figure \ref{fig:a1_m50_beta1}. It can be seen that no stability problems arise, since the real part of all eigenvalues are non positive.

\subsection*{Strong form, collocation approach with SGLR  nodes and scaled Laguerre polynomials}

\begin{figure}[htb]
\subfloat[]{\includegraphics[scale=0.45]{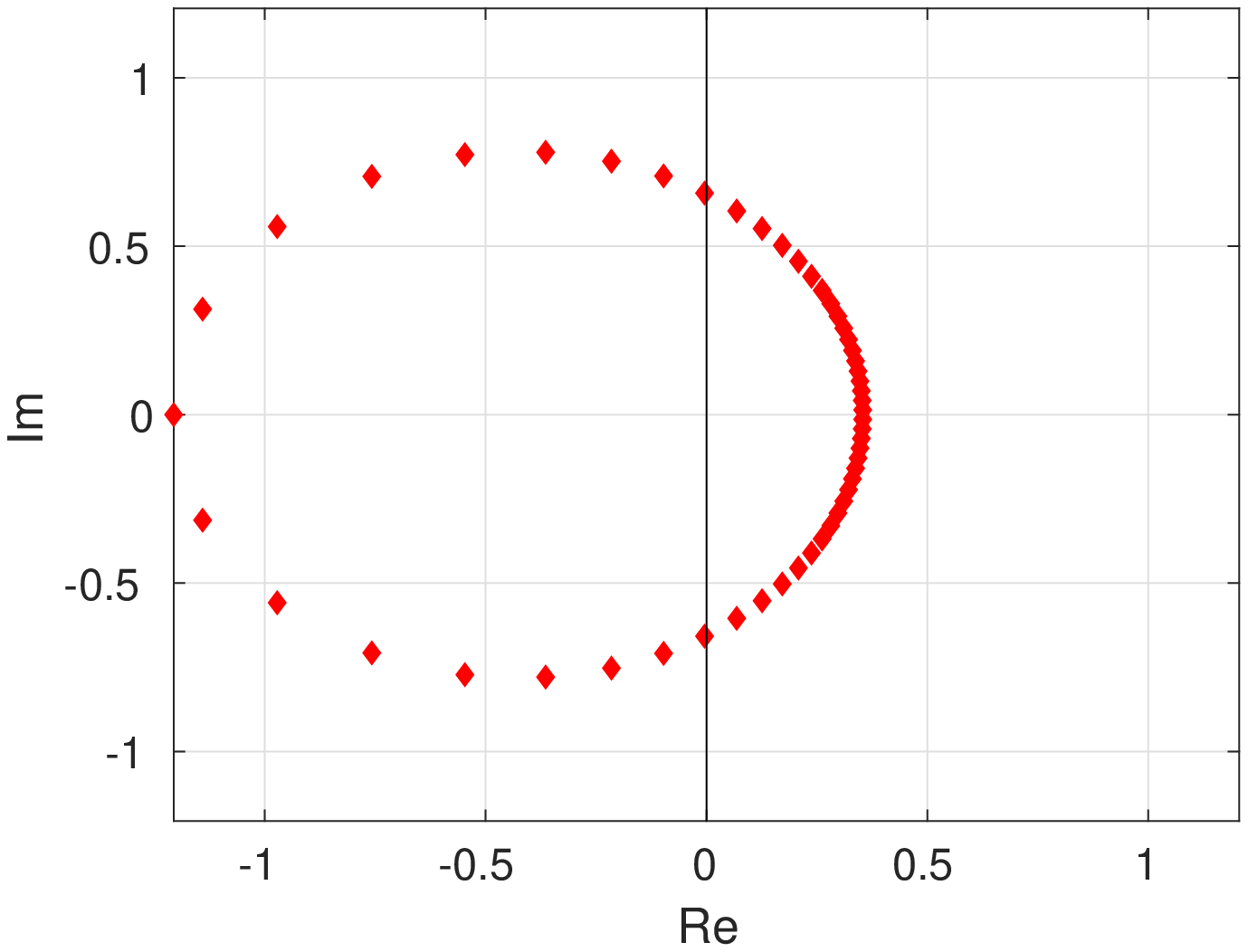}}\quad
\subfloat[]{\includegraphics[scale=0.45]{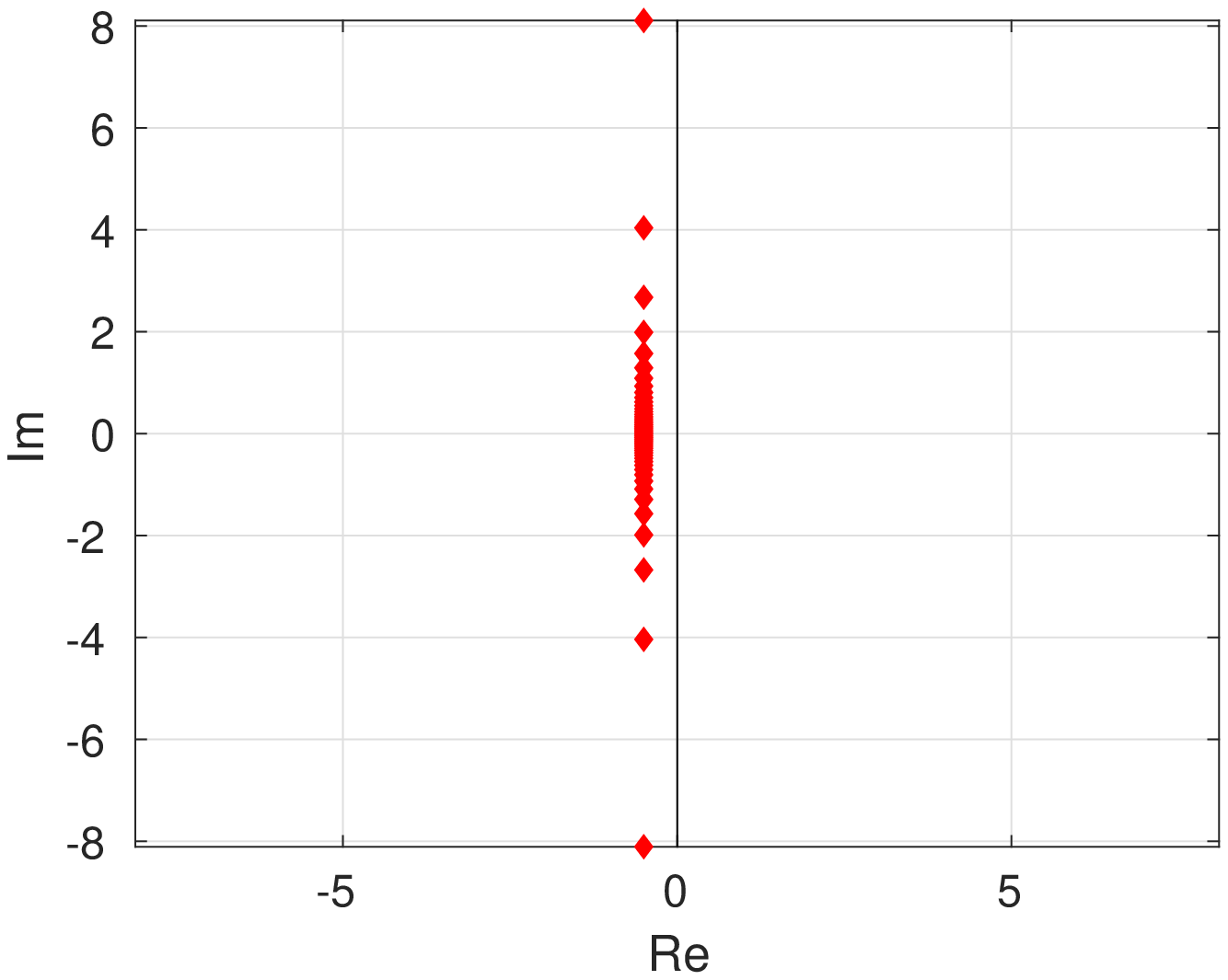}}
\caption{Eigenvalues in the a1 case, scaled Laguerre polynomials, for parameter values  $\beta=1$, $M=50$
in the a) outflow  case $u=-1$ and b) inflow  case $u=1$. }
\label{fig:a1_m50_beta1_pol}
\end{figure}

If scaled Laguerre polynomials instead of Laguerre functions are used (see \cite{shen:2001, shen:2009} for definitions),
one gets the discretized equations \eqref{coll_sglr_in} and \eqref{coll_sglr_out} with, instead of $\hat h^\beta_j$, the Lagrange interpolation
functions $h^\beta_j$ associated with the nodes and weights obtained considering the scaled Laguerre polynomials basis,
and $\mathbf{D}_\beta$ now denoting the SGLR differentiation matrix relative to that basis. 
 Examples of eigenvalue
plots for $\mathbf{A}$ in this case are shown in figure \ref{fig:a1_m50_beta1_pol}. It can be seen
that the inflow case does not display any stability
issue, while in the outflow case matrix has some eigenvalues with positive real part.
\subsection*{Weak form, nodal approach with scaled Laguerre functions}

For the case of weak form, nodal discretizations based on scaled Laguerre functions,
 we take again as test and basis functions
the Lagrange interpolation functions associated with the SGL or SGLR nodes. 
We then use the fact that
\begin{equation}
q(z)\approx \sum_{j=0}^M q_j \hat h^{\beta}_j(z). \label{nodal_spec_exp}
\end{equation}
We consider the   SGL nodes first.
Using the fact that the basis functions are also a Lagrange basis and employing the
corresponding numerical integration rule,  in the inflow case we get from 
\eqref{weak_form_ibp_in} 
\begin{equation}
q_i^{\prime}\omega_i^{\beta}= uq_L\hat h^{\beta}_i(0) +u\sum_{j=0}^M q_j (\hat h^{\beta}_j)^{\prime}(z_i^{\beta})
\omega_j^{\beta} \quad i=0,\dots,M
\label{weak_nodal_sgl_in}
\end{equation}
and in the outflow case from \eqref{weak_form_ibp_out} 
\begin{equation}
q_i^{\prime}\omega_i^{\beta}= u\sum_{j=0}^M q_j \hat h^{\beta}_j(0)\hat h^{\beta}_i(0) 
+u\sum_{j=0}^M q_j (\hat h^{\beta}_j)^{\prime}(z_i^{\beta})
\omega_j^{\beta} \quad i=0,\dots,M.
\label{weak_nodal_sgl_out}
\end{equation}
In matrix notation, equation \eqref{weak_nodal_sgl_in} yields in the inflow case 
an ODE system like equation \eqref{disc_advection} with
\begin{equation}
\mathbf{A}=u\mathbf{\Omega}^{-1}\mathbf{D}_{\beta}\mathbf{\Omega},
\end{equation}
where $\mathbf{\Omega}$ is the diagonal matrix with the weights on the main diagonal,
and
\begin{equation}
\mathbf{g} = uq_L\mathbf{\Omega}^{-1}\mathbf{h},
\end{equation}
where we have set $\mathbf{h}=[\hat h^{\beta}_0(0),\dots,\hat h^{\beta}_M(0)] ^T$.
For the outflow case,  equation \eqref{weak_nodal_sgl_out} yields instead 
$\mathbf{g}=\mathbf{0}$ and
\begin{equation}
\mathbf{A}=u \mathbf{\Omega}^{-1}\mathbf{H}+u\mathbf{\Omega}^{-1}\mathbf{D}_{\beta}\mathbf{\Omega},
\end{equation}
where $ \mathbf{H} = \mathbf{h}\mathbf{h}^T$.  

\begin{figure}[htb]
\subfloat[]{\includegraphics[scale=0.45]{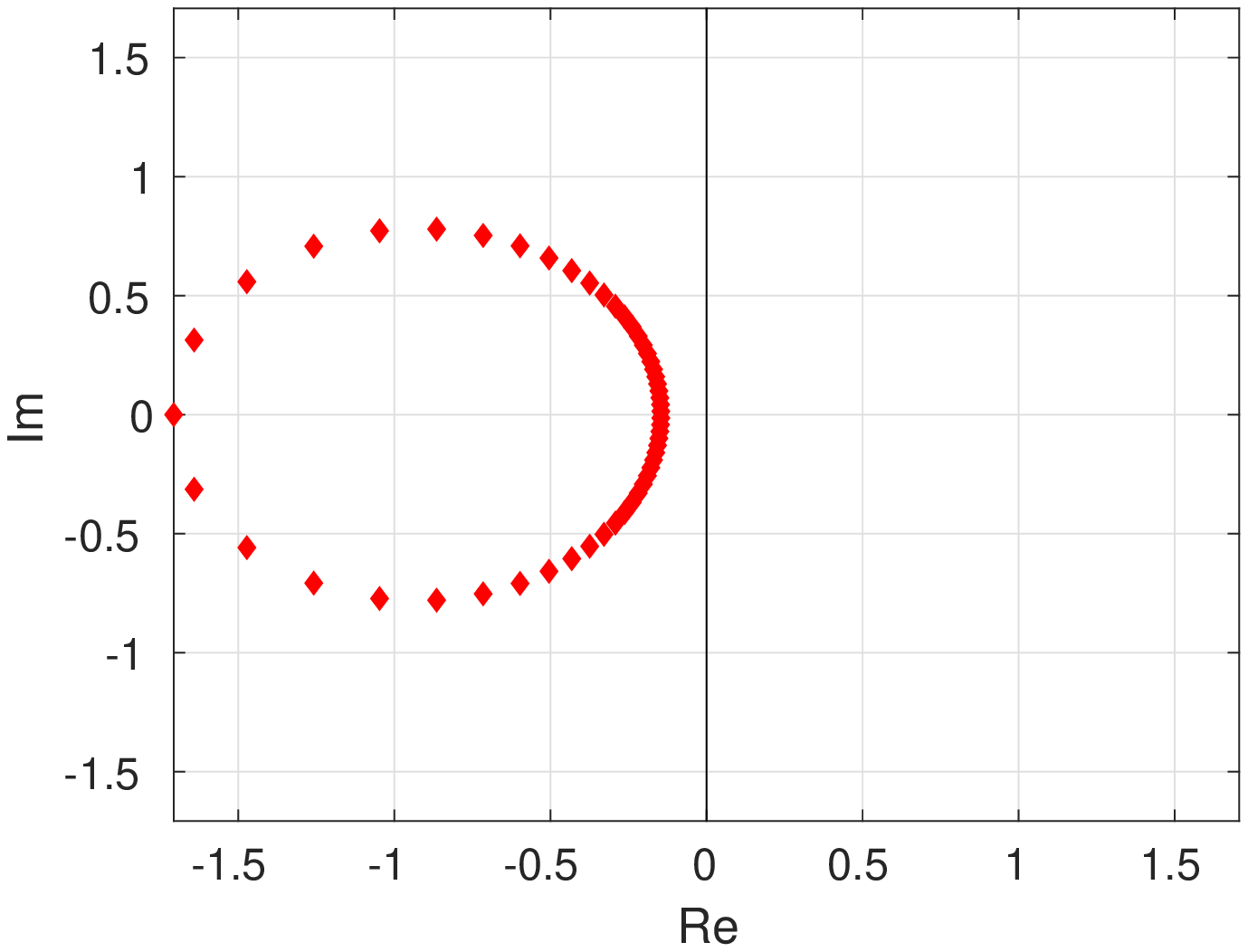}}\quad
\subfloat[]{\includegraphics[scale=0.45]{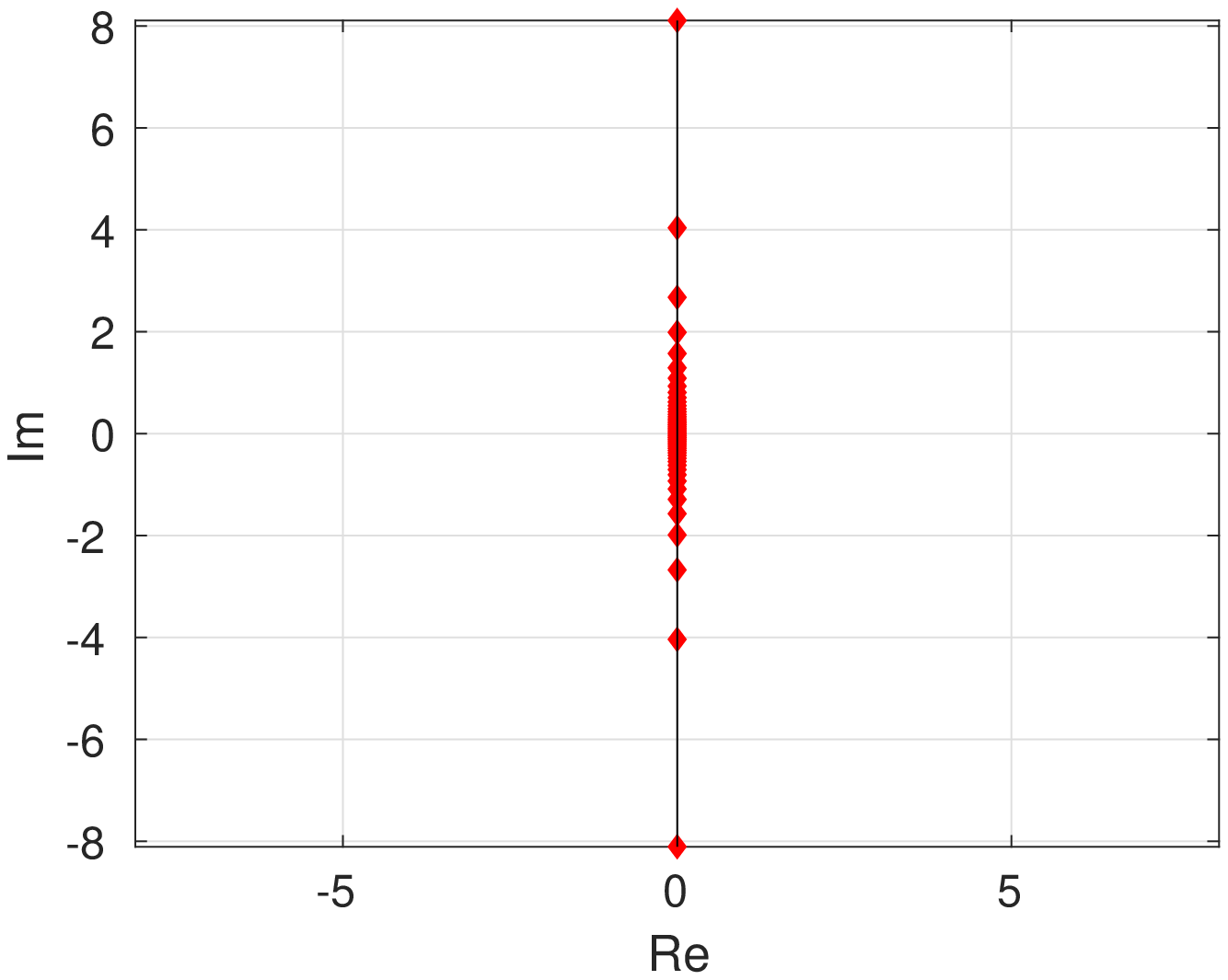}}
\caption{Eigenvalues in the a2n case with GLR nodes, for parameter values  $\beta=1, $ $M=50$
in the a) outflow  case $u=-1$ and b) inflow  case $u=1$. }
\label{fig:a2n_glr_m50_beta1}
\end{figure}

\begin{figure}[htb]
\centering
\includegraphics[scale=0.45]{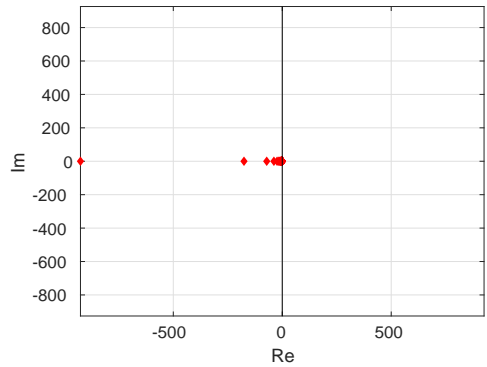}  
\caption{Eigenvalues in the a2n case with GL nodes, for parameter values  $\beta=1, $ $M=50$
in the  inflow  case $u=1$. }
\label{fig:a2n_gl_m50_beta1}
\end{figure}

 For the case of SGLR nodes, one  has in the inflow case, 
\begin{equation}
q_i^{\prime}\omega_i^{\beta} =   uq_L (\hat h^{\beta}_0)^{\prime}(z_i^{\beta}) \omega_0^{\beta}
+u\sum_{j=1}^M q_j (\hat h^{\beta}_j)^{\prime}(z_i^{\beta})
\omega_j^{\beta}  \quad i=1,\dots,M.
\label{weak_nodal_sglr_in}
\end{equation}
where we have used   the fact that for the SGLR Lagrangian basis one has $\hat h^{\beta}_i(0)=0$ for $ i=1,\dots,M $
while $\hat h^{\beta}_0(0)=1$. In the outflow case, one has instead
\begin{equation}
q_i^{\prime}\omega_i^{\beta}=   u\sum_{j=0}^M q_j (\hat h^{\beta}_j)^{\prime}(z_i^{\beta})\omega_j^{\beta}, \quad i=1,\dots,M
\label{weak_nodal_sglr_out}
\end{equation}
while for $i=0 $ one obtains:
\begin{equation}
q_0^{\prime}\omega_0^{\beta}= uq_0 +u\sum_{j=0}^M q_j (\hat h^{\beta}_j)^{\prime}(z_0^{\beta})\omega_j^{\beta}.
\label{weak_nodal_sglr_out0}
\end{equation}
In matrix notation, equation \eqref{weak_nodal_sglr_in} yields in the inflow case an ODE system like equation \eqref{disc_advection} with $\mathbf{q}=[q_1,\dots,q_M]^T$,  
\begin{equation}
\mathbf{g}=\omega_0^{\beta}uq_L\left[ \frac{(\hat h^{\beta}_0)^{\prime}(z_1^{\beta})}{\omega_1^{\beta}},\dots, \frac{(\hat h^{\beta}_0)^{\prime}(z_M^{\beta})}{\omega_M^{\beta}}\right] ^T
\end{equation}
and $\mathbf{A}=u\mathbf{\Omega}^{-1}_M(\mathbf{D}_{\beta})_M \mathbf{\Omega}_M$, where again $\mathbf{\Omega}_M $ denotes the $M\times M$ diagonal matrix obtained omitting
from $\mathbf{\Omega}$ the first row and the first column. In the outflow case instead one has $\mathbf{q}=[q_0,\dots,q_M]^T, $ $\mathbf{g}=\mathbf{0} $ and 
\begin{equation}
\mathbf{A}=\frac{u}{\omega^{\beta}_0}\mathbf{e}_1\mathbf{e}_1^T+u\mathbf{\Omega}^{-1}\mathbf{D}_{\beta}\mathbf{\Omega}.
\end{equation}
Examples of eigenvalue plots for $\mathbf{A} $ in the a2n case with GLR nodes are shown in figure \ref{fig:a2n_glr_m50_beta1}.  
It can be seen that the eigenvalue structure is entirely analogous to that of the collocation case a1.
This contrasts with the analogous case discretized  with GL nodes. Indeed, it can be seen in figure \ref{fig:a2n_gl_m50_beta1} that eigenvalues are much more spread out,
with large negative real parts that identify a potentially stiff problem. The situation is even worse for the outflow case (not shown), for which the largest eigenvalue
(in absolute value) appears to be of the order $10^{52}$.
\subsection*{Weak form, nodal approach  with scaled Laguerre polynomials}
Using Laguerre polynomials, the expansion \eqref{nodal_spec_exp} reads instead:
\begin{equation}
q(z)\approx \sum_{j=0}^M q_jh^{\beta}_j(z). \label{nodal_spec_exp_pol}
\end{equation}
Considering SGL nodes first, and replacing the expression in \eqref{weak_form_poly_ibp_in} we find, in the inflow case:
\begin{multline}
\frac{d}{dt}\int_0^{+\infty}h^\beta_i(z)\sum_{j=0}^Mq_jh^\beta_j(z)\omega(z)\,dz +\beta u\int_0^{+\infty}h^\beta_i(z)\sum_{j=0}^Mq_jh^\beta_j(z)\omega(z) \,dz \\
 -uq_Lh^\beta_i(0)\omega(0) -u\int_0^{+\infty} \left(h^\beta_i\right)'(z)\sum_{j=0}^Mq_jh^\beta_j(z)\omega(z) \,dz =0
\end{multline}
and in the outflow case, using \eqref{weak_form_poly_ibp_out}:
\begin{multline}
\frac{d}{dt}\int_0^{+\infty}h^\beta_i(z)\sum_{j=0}^Mq_jh^\beta_j(z)\omega(z)\,dz  +\beta u\int_0^{+\infty}h^\beta_i(z)\sum_{j=0}^Mq_jh^\beta_j(z)\omega(z) \,dz \\
- u\sum_{j=0}^Mq_jh_j^\beta(0)h^\beta_i(0)\omega(0) - u\int_0^{+\infty} \left(h^\beta_i\right)'(z)\sum_{j=0}^Mq_jh^\beta_j(z)\omega(z) \,dz =0
\end{multline}
Discretizing the integrals with quadrature formulas based on Laguerre polynomials we get, for $i=0,\ldots,M$, in the inflow case:
\begin{multline}
\frac{d}{dt}\sum_{k=0}^Mh^\beta_i(z_k^\beta)\sum_{j=0}^Mq_jh^\beta_j(z_k^\beta)\omega_k^\beta + \beta u\sum_{k=0}^Mh^\beta_i(z_k^\beta)\sum_{j=0}^Mq_jh^\beta_j(z_k^\beta)\omega_k^\beta \\
 - uq_Lh^\beta_i(0)\omega(0) -u\sum_{k=0}^M\left(h^\beta_i\right)'(z_k^\beta)\sum_{j=0}^Mq_jh^\beta_j(z_k^\beta)\omega_k^\beta =0
\end{multline}
and a similar expression for the outflow case upon replacement of $q_L$ with $\sum_{j=0}^Mq_jh_j^\beta(0)$.
With further simplifications and using also the fact that $\omega(0)=1, $ we get for $i=0,\ldots,M$ in the inflow case:
\begin{equation}
q_i'\omega_i^\beta = uq_Lh^\beta_i(0) - \beta u q_i\omega_i^\beta + u\sum_{j=0}^Mq_j(h_j^\beta)'(z_i^\beta)\omega_j^\beta\\
\end{equation}
and in the outflow case:
\begin{equation}
q_i'\omega_i^\beta = u\sum_{j=0}^Mq_jh_j^\beta(0)h^\beta_i(0) - \beta u q_i\omega^\beta_i + u\sum_{j=0}^Mq_j(h_j^\beta)'(z_i^\beta)\omega_j^\beta
\end{equation}
In matrix form, in the inflow case we have:
\begin{equation}
\mathbf{A} = u\mathbf{\Omega}^{-1}\mathbf{D}_\beta\mathbf{\Omega} 
- \beta u \mathbf{I},\quad\mathbf{g}=uq_L\mathbf{\Omega}^{-1}\mathbf{h}, \quad\mathbf{h}=[h_0(0),\ldots,h_M(0)]^T
\end{equation}
while in the outflow case $\mathbf{g}=0$ and
\begin{equation}
\mathbf{A} = u\mathbf{\Omega}^{-1}\mathbf{H}  +u\mathbf{\Omega}^{-1}\mathbf{D}_\beta\mathbf{\Omega} - \beta u \mathbf{I}.
\end{equation}
The spectrum of this matrix behaves in a manner that is analogous to the weak form, Laguerre function case
shown in figure \ref{fig:a2n_gl_m50_beta1}. Also in this case, the eigenvalues in the outflow case become
extremely large.

\begin{figure}[htb]
\subfloat[]{\includegraphics[scale=0.45]{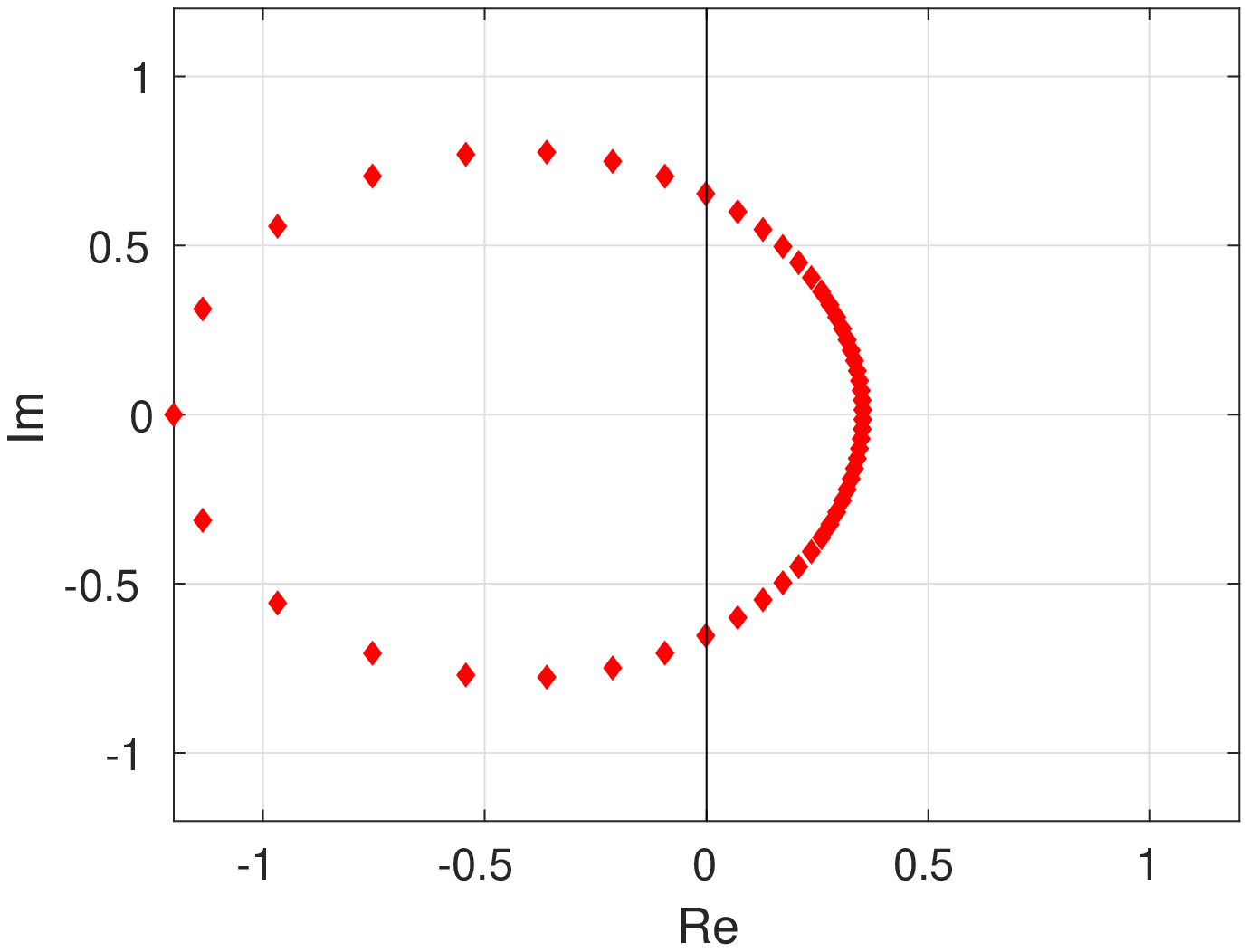}}\quad
\subfloat[]{\includegraphics[scale=0.45]{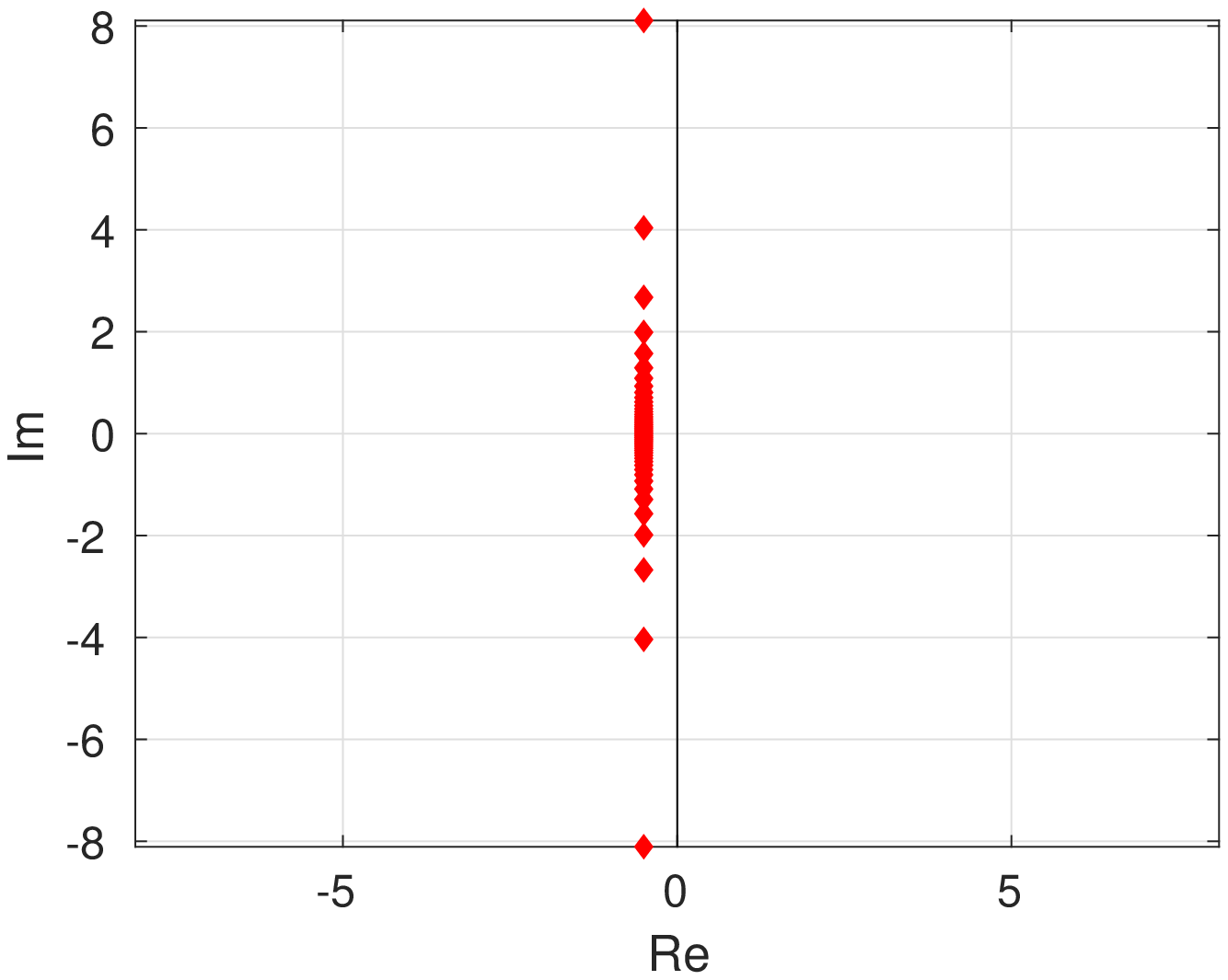}}
\caption{Eigenvalues in the a2n case with GLR nodes, scaled Laguerre polynomials, for parameter values  $\beta=1, $ $M=50$
in the a) outflow  case $u=-1$ and b) inflow  case $u=1$. }
\label{fig:a2n_glr_m50_beta1_pol}
\end{figure}

In the case of SGLR nodes, replacing in \eqref{weak_form_poly_ibp_in} we find in the inflow case, for $i=1,\ldots,M$:
\begin{equation}
q_i'\omega^\beta_i =  uq_L\left(h^\beta_0\right)'(0)\omega_0^\beta-\beta u q_i\omega_i + u\sum_{j=1}^Mq_j\left(h_j^\beta\right)'(z_i^\beta)\omega_j^\beta
\end{equation}
Replacing in \eqref{weak_form_poly_ibp_out} we find in the outflow case, for $i=1,\ldots,M$
\begin{equation}
q_i'\omega^\beta_i = -\beta u q_i\omega_i + u\sum_{j=0}^Mq_j\left(h_j^\beta\right)'(z_i^\beta)\omega_j^\beta
\end{equation}
while for $i=0$ we get:
\begin{equation}
q_0'\omega^\beta_0 = uq_0 -\beta u q_0\omega_0 + u\sum_{j=0}^Mq_j\left(h_j^\beta\right)'(0)\omega_j^\beta.
\end{equation}
In matrix form, in the inflow case we have, denoting  $\mathbf{q}=[q_1,\dots,q_M]^T$, 
\begin{align}
\mathbf{A}&= u\mathbf{\Omega}_M^{-1}\left(\mathbf{D}_\beta\right)_M\mathbf{\Omega}_M
-\beta u\mathbf{I}_M,\nonumber\\
\mathbf{g}&= \omega_0^\beta uq_L\left[\frac{\left(h^\beta_0\right)'\left(0\right)}{\omega_1^\beta},\ldots,
\frac{\left(h^\beta_0\right)'\left(0\right)}{\omega_M^\beta}\right]^T
\end{align}
In the outflow case, $\mathbf{q}=[q_0,\dots,q_M]^T$, $\mathbf{g}=0$ and 
\begin{equation}
\mathbf{A} = \frac{u}{\omega_0^\beta}\mathbf{e_1}\mathbf{e_1}^T+u\mathbf{\Omega}^{-1}\mathbf{D}_\beta\mathbf{\Omega}-\beta u \mathbf{I}
\end{equation} 
The spectrum of the matrix $\mathbf{A}$ is depicted in figure \ref{fig:a2n_glr_m50_beta1_pol}, the findings are equivalent to the collocation
discretization of the model in strong form.

\subsection*{Weak form, modal approach with scaled Laguerre  functions}

In case a modal discretization based on the scaled Laguerre  functions is sought,
the solution will instead be represented as
\begin{equation}
q(z)\approx \sum_{j=0}^M  q_j\hat{\mathscr L}^\beta_j(z),
\end{equation}
where
\begin{equation}
  q_j=\beta \int_0^{+\infty}  q(z) \hat{\mathscr L}^\beta_j(z)\,dz. 
\end{equation}
This implies that, substituting in \eqref{weak_form_ibp_in}  and choosing $\hat{\mathscr L}^\beta_i $
as test function we obtain, in the inflow case,
\begin{equation}
\frac1\beta  q_i^{\prime} = uq_L 
+u\sum_{j=0}^M   q_j \int_0^{+\infty} \frac{\partial }{\partial z}\hat{\mathscr L}^\beta_i(z) \hat{\mathscr L}^\beta_j(z)\,dz \quad i=0,\dots,M,
\label{weak_modal_fun_in0}
\end{equation}
where we have used the fact that $\hat{\mathscr L}^\beta_i (0)=1$. 
One can then consider that 
\begin{equation}
\label{der_laguerre}
\frac{\partial }{\partial z}\hat{\mathscr L}^\beta_i(z) =-\frac{\beta}2 \hat{\mathscr L}^\beta_i(z)-\beta
\sum_{k=0}^{i-1}\hat{\mathscr L}^\beta_k(z) 
\end{equation}
to obtain
\begin{align}
\sum_{j=0}^M   q_j \int_0^{+\infty}   \frac{\partial }{\partial z}\hat{\mathscr L}^\beta_i(z) \hat{\mathscr L}^\beta_j(z)\,dz
= &-\frac{\beta}2\sum_{j=0}^M   q_j \int_0^{+\infty} \hat{\mathscr L}^\beta_i(z)\hat{\mathscr L}^\beta_j(z)\,dz
\nonumber
\\
 & -\beta \sum_{j=0}^M q_j \int_0^{+\infty}\hat{\mathscr L}^\beta_j(z) \sum_{k=0}^{i-1}\hat{\mathscr L}^\beta_k(z)\,dz \nonumber \\
= &- \frac{q_i}2 - \sum_{j=0}^{i-1}  q_j,
\label{weak_modal_fun_in1}
\end{align}
so that in the end we obtain 
\begin{equation}
  q_i^{\prime} = \beta uq_L    
-\beta u\frac{  q_i}2  -\beta u\sum_{j=0}^{i-1}   q_j \qquad i=0,\dots,M.
\label{weak_modal_fun_in2}
\end{equation}
Going through the same steps after substitution in \eqref{weak_form_ibp_out}, we obtain for $i=0,\dots,M$
\begin{equation}
\frac1\beta  q_i^{\prime} = u \sum_{j=0}^M  q_j
+u\sum_{j=0}^M   q_j \int_0^{+\infty} \frac{\partial }{\partial z}\hat{\mathscr L}^\beta_i(z) \hat{\mathscr L}^\beta_j(z)\,dz \label{weak_modal_fun_out0}
\end{equation}
that simplifies to 
\begin{eqnarray}
   q_i^{\prime} &=& \beta u\sum_{j=0}^M  q_j   
-\beta u\frac{  q_i}2 -\beta u\sum_{j=0}^{i-1}   q_j \nonumber \\
&=&\beta u \sum_{j=i+1}^M  q_j +\beta u\frac{ q_i}2
\quad i=0,\dots,M.
\label{weak_modal_fun_out1}
\end{eqnarray}
In matrix notation, equation \eqref{weak_modal_fun_in2} yields for the inflow case 
an ODE system like equation \eqref{disc_advection} with
$\mathbf{q}=[q_0,\dots,q_M]^T$, $\mathbf{g}=uq_L[1,\dots,1]^T$ and $\mathbf{A}=-\beta u\mathbf{L}$,
where $\mathbf{L}$ denotes a lower triangular $(M+1)\times(M+1)$ with values equal to one half on the main diagonal and one below.
In the outflow case instead one has $\mathbf{q}=[q_0,\dots,q_M]^T$, $\mathbf{g}=\mathbf{0}$ and $\mathbf{A}=\beta u\mathbf{L}^{T}$.
In both cases, the eigenvalues of the matrix are all equal to $-\beta |u|/2$, so that no stability problems arise.

\subsection*{Weak form, modal approach with scaled Laguerre  polynomials}
If a modal discretization based on the scaled Laguerre polynomials is sought, one can go through the similar steps as in the previous Sections assuming
\begin{equation}
q(z)\approx \sum_{j=0}^M  q_j {\mathscr L}^\beta_j(z),
\end{equation}
and taking into account that the recurrence relationship for the derivatives
now reads:
\begin{equation}
 \frac{\partial }{\partial z}{\mathscr L}^\beta_i(z) = 
 -\beta \sum_{k=0}^{i-1}{\mathscr L}^\beta_k(z)  \quad
 {\rm for} \ i\geq 1 \quad  {\rm and} \quad \frac{\partial }{\partial z}{\mathscr L}^\beta_0(z)=0.
\end{equation}
It follows that, for $i\geq 1$,
\begin{align}
&\sum_{j=0}^M  q_j \int_0^{+\infty} \frac{\partial }{\partial z}{\mathscr L}^\beta_i(z){\mathscr L}^\beta_j(z)\omega(z)\,dz=
\nonumber
\\
& -\beta \sum_{j=0}^M   q_j \int_0^{+\infty}{\mathscr L}^\beta_j(z) \sum_{k=0}^{i-1}{\mathscr L}^\beta_k(z)\omega(z)\,dz =  -\sum_{j=0}^{i-1}   q_j,
\label{weak_modal_pol_in1}
\end{align}
so that in the end, after substitution in \eqref{weak_form_poly_ibp_in}, we obtain 
\begin{gather}
\frac1\beta  q_0^{\prime} = uq_L -  u q_0   
\label{weak_modal_pol_in20}\\
 \frac1\beta  q_i^{\prime} = uq_L - u  q_i   
 - u\sum_{j=0}^{i-1} q_j = uq_L - u\sum_{j=0}^i   q_j \qquad i=1,\dots,M.
\label{weak_modal_pol_in2}
\end{gather}
Going through the same steps after substitution in \eqref{weak_form_poly_ibp_out}, we obtain for $i=0,\dots,M$
\begin{equation}
\frac1\beta q_i^{\prime} = u\sum_{j=0}^M  q_j - u  q_i
+u\sum_{j=0}^M   q_j \int_0^{+\infty} \frac{\partial }{\partial z}{\mathscr L}^\beta_i(z) {\mathscr L}^\beta_j(z)\,dz   \label{weak_modal_pol_out0}
\end{equation}
that simplifies to 
\begin{gather}
\frac1\beta   q_0^{\prime} =   u    \sum_{j=0}^M  q_j   
   - u  q_0     \label{weak_modal_pol_out10}\\
\frac1\beta  q_i^{\prime} =   u    \sum_{j=0}^M q_j   
   - u q_i -  u\sum_{j=0}^{i-1}  q_j
= u \sum_{j=i+1}^M  q_j    
   \qquad i=1,\dots,M.
   \label{weak_modal_pol_out1}
\end{gather}
In matrix notation, equation \eqref{weak_modal_pol_in2} yields for the inflow case an ODE system like equation \eqref{disc_advection} with
$\mathbf{q}=[q_0,\dots, q_M]^T$, $\mathbf{g}=uq_L[1,\dots,1]^T$ and $\mathbf{A}=-\beta u\mathbf{L}$, where $\mathbf{L}$ denotes a lower triangular  $(M+1)\times (M+1)$ matrix 
 with values equal to one on and below the main diagonal. In the outflow case instead one has $\mathbf{q}=[q_0,\dots,q_M]^T, $ $\mathbf{g}=\mathbf{0}$ and 
 $\mathbf{A}=\beta u\mathbf{U}$, where   $\mathbf{U}$ denotes an upper triangular  $(M+1)\times (M+1)$ matrix with values equal to zero on the main diagonal and one above the main diagonal.
In the inflow case, the eigenvalues of the matrix are all equal to $ -\beta u$, so that no stability problems arise. In the outflow case, the matrix 
has all eigenvalues equal to zero, so that again no stability problems arise.

\section{A scaled Laguerre functions discretization of hyperbolic systems}
\label{modal_sys}

In this Section, we present a  discretization of hyperbolic systems on semi-infinite domains.
Based on the analysis in Section \ref{analysis}, we opt for a weak form discretization with a
 modal basis of Laguerre scaled functions. In order to account for cases of practical interest
 in relevant applications, we consider the system

\begin{equation}
\displaystyle \frac{\partial \mathbf{q}}{\partial t} 
+    \mathbf{A}(\mathbf{q},z)  \frac{\partial \mathbf{q}}{\partial z}   
 = \mathbf{B}(\mathbf{q},z) \mathbf{q}.
\label{lin1d}
\end{equation}
The system is hyperbolic if, for any $(\mathbf{q},z), $ the matrix
  $ \mathbf{A}(\mathbf{q},z)   \in \mathbb{R}^d \times \mathbb{R}^d$
is diagonalizable as  $ \mathbf{A} = \mathbf{V} \boldsymbol{\Lambda} \mathbf{V}^{-1}$, and has real eigenvalues.
The matrix $\mathbf{B}$ represents a zero order forcing term, such as for example the Coriolis force or buoyancy
effects in the propagation of internal gravity waves.
We first consider equations \eqref{lin1d} on the whole domain  $ {\mathbb R}^+$, written componentwise as
\begin{equation}
 \frac{\partial q_k}{\partial t} 
+   \sum_{l=1}^d a_{kl}   \frac{\partial q_l}{\partial z}   
 = \sum_{l=1}^d b_{kl} q_l, \ \ \ k=1,\dots,d
\label{lin1d_comp}
\end{equation}
where $a_{kl},  b_{kl} $ are the entries of  $ \mathbf{A}, $  $ \mathbf{B}$, respectively.
Notice that the dependency on  $\mathbf{q}, z, $ has been omitted to simplify the notation.
 We then define as customary the matrices
$ \mathbf{A}^{\pm}=\mathbf{V} \boldsymbol{\Lambda}^{\pm} \mathbf{V}^{-1}$, where $\lambda_i^{+}= \max\{\lambda_i,0\}, $ $\lambda_i^{-}=  \min\{\lambda_i,0\}. $
We denote by $g_{k}(t) $ the value assigned to each component of the system by the Dirichlet boundary conditions and by $\hat q_{k,j}, j=0,\dots,M $ the $j-$th modal coefficient of $q_k$,
so that each component will be represented as
\begin{equation}
\label{modalexp}
q_k(z,t)\approx \sum_{j=0}^M q_{k,j}(t)\hat{\mathscr L}^\beta_j(z),
\end{equation}
where
\begin{equation}
  q_{k,j}=\beta \int_0^{+\infty}  q_k(z) \hat{\mathscr L}^\beta_j(z)\,dz. 
\end{equation}
Integrating \eqref{lin1d_comp} against a test function $\hat{\mathscr L}^\beta_i,$
we have again for $ k=1,\dots,d$ 

\begin{eqnarray}
 \int_0^{+\infty}   \frac{\partial q_k}{\partial t} (z,t)\hat{\mathscr L}^\beta_i (z) \,dz
&=&  -  \sum_{l=1}^d   \int_0^{+\infty}  \hat{\mathscr L}^\beta_i(z) a_{kl}   \frac{\partial q_l}{\partial z}(z,t)   \,dz
\nonumber \\
 &+&  \sum_{l=1}^d  \int_0^{+\infty}  \hat{\mathscr L}^\beta_i (z) b_{kl} q_l(z,t), \,dz. 
\nonumber 
\end{eqnarray}
Since, performing  integration by parts, one has
\begin{eqnarray}
&& \sum_{l=1}^d   \int_0^{+\infty}  \hat{\mathscr L}^\beta_i(z) a_{kl}  \frac{\partial q_l}{\partial z}(z,t)  \,dz
\nonumber \\
 && = - \hat{\mathscr L}^\beta_i(0)  \sum_{l=1}^d a^+_{kl}g_l(t) 
 - \hat{\mathscr L}^\beta_i(0)  \sum_{l=1}^d a^-_{kl}q_l(0)\nonumber \\
&& -  \sum_{l=1}^d   \int_0^{+\infty}  q_l (z,t)a_{kl}  \frac{\partial \hat{\mathscr L}^\beta_i}{\partial z}  (z)  \,dz \nonumber\\
&& -  \sum_{l=1}^d   \int_0^{+\infty}  q_l (z,t) \hat{\mathscr L}^\beta_i(z) \frac{\partial  a_{kl} }{\partial z}  (z)  \,dz \nonumber\\
 &&+ \sum_{l=1}^d  \int_0^{+\infty}  \hat{\mathscr L}^\beta_i (z) b_{kl} q_l(z,t), \,dz. 
\label{lin1d_comp_int2}
\end{eqnarray}
using the fact that $\hat{\mathscr L}^\beta_i(0)=1$ it follows

\begin{eqnarray}
\frac{1}{\beta}\frac{d}{dt}  q_{k,i}(t)&=&\sum_{l=1}^d a^{+}_{kl} g_{l}(t)+
\sum_{l=1}^d a^{-}_{kl} \sum_{j=0}^M   q_{l,j}(t)   \nonumber\\
&& +  \sum_{l=1}^d   \int_0^{+\infty}  q_l (z,t)a_{kl}  \frac{\partial \hat{\mathscr L}^\beta_i}{\partial z}  (z)  \,dz \nonumber\\
&& + \sum_{l=1}^d   \int_0^{+\infty}  q_l (z,t) \hat{\mathscr L}^\beta_i(z) \frac{\partial  a_{kl} }{\partial z}  (z)  \,dz \nonumber\\
 &&+ \sum_{l=1}^d  \int_0^{+\infty}  \hat{\mathscr L}^\beta_i (z) b_{kl} q_l(z,t) \,dz. 
\label{modal_sys_1}
\end{eqnarray}
Using formulae \eqref{der_laguerre}  and \eqref{modalexp} one obtains then
\begin{eqnarray}
\frac{1}{\beta}\frac{d}{dt} q_{k,i}(t)&=&\sum_{l=1}^d a^{+}_{kl} g_{l}(t)
+\sum_{l=1}^d a^{-}_{kl} \sum_{j=0}^M  q_{l,j}(t)   \nonumber\\
&& - \frac{\beta}2 \sum_{l=1}^d  
\sum_{j=0}^M  q_{l,j}(t) \int_0^{+\infty} a_{kl} \hat{\mathscr L}^\beta_i (z) \hat{\mathscr L}^\beta_j (z) \,dz \nonumber\\
&& - \beta \sum_{l=1}^d  
\sum_{j=0}^M   q_{l,j}(t) \sum_{s=0}^{i-1} \int_0^{+\infty} a_{kl} \hat{\mathscr L}^\beta_s (z) \hat{\mathscr L}^\beta_j (z) \,dz \nonumber\\
&& +  \sum_{l=1}^d  \sum_{j=0}^M  q_{l,j}(t) \int_0^{+\infty}  \hat{\mathscr L}^\beta_i(z) \hat{\mathscr L}^\beta_j (z) \frac{\partial  a_{kl} }{\partial z}  (z)  \,dz \nonumber\\
 &&+ \sum_{l=1}^d  \sum_{j=0}^M  q_{l,j}(t)\int_0^{+\infty}  b_{kl} \hat{\mathscr L}^\beta_i (z) \hat{\mathscr L}^\beta_j (z) \,dz.
 \label{modal_sys_2}
\end{eqnarray}
In the linear, constant coefficient case, the integrals in \eqref{modal_sys_2} yield expressions similar to those
derived in Section \ref{analysis}, while in the variable coefficient or nonlinear case
a full space discretization is  obtained by applying the Gauss - Laguerre - Radau quadrature formulae.

%
%

\section{Coupling of the scaled Laguerre discretization
  with a DG discretization on the finite domain}
\label{coupled}

We now split the domain ${\mathbb R}^+ $ as ${\mathbb R}^+=[0,L]\cup[L,+\infty)$ in order to introduce two different discretizations on the finite and semi-infinite parts of the
domain. On $[L,+\infty), $ system  \eqref{lin1d_comp} is discretized by \eqref{modal_sys_2}, where
an appropriate shift is performed in the independent variable. 

On the finite domain, a DG discretization is employed. For definiteness, we have applied here
a  modal DG approach.
More specifically, a computational mesh is introduced in the $[0,L] $ interval by defining a set of $N $ non overlapping elements $K_m$ of size $\Delta z_m$, such that
$[0,L]=\bigcup_{m=1}^{N} K_m$.
The center of the generic element $K_m$ is denoted by $z_m$, while $z_{m\pm1/2} $ denote its boundary points. Each element  $K_m$ can be seen as the image of the master element
$\hat K = [-1,1]  $ via the affine local map $z = f_{m}(\xi)= \xi\Delta z_m/2+z_{m}$, where $z\in K_m $ and $\xi\in \hat K$.
For each non-negative integer $p$, we then denote by $\mathbb{P}_p$ the set of all polynomials of degree less or equal to $p$ on $\hat K$. We will also define
$\mathbb{P}_p(K_m)= \left\{w : w= v\circ {F}^{-1}_{m}, \quad  v\in \mathbb{P}_p \right\}$.
For each polynomial degree $p$, the   discontinuous finite element spaces are defined as follows
\begin{equation}
V^p_{h}=\left\{ v \in L^2( [0,L] ) :
v|_{K_m} \in \mathbb{P}_{p}(K_m)\quad  m=1,\dots,N  \right\},
\label{fespace_disc}
\end{equation}
For each element $K_m,\,m=1,\dots,N$  we then denote by $\phi^{m}_j(z), j=0,\dots, p$  a basis of $ \mathbb{P}_{p}(K_m)$.  
For discontinuous finite elements discretizations, we will consider instead the orthogonal basis based on Legendre polynomials.
More specifically, for $\xi \in \hat{K}$, define the Legendre polynomial recursively by the following recurrence relation:
\begin{equation}
\begin{array}{ll}
\displaystyle L_{k+1} = \frac{2k+1}{k+1} \xi L_k(\xi) - \frac{k}{k+1} L_{k-1}(\xi), & k = 1,2,\ldots \bigskip \\
L_0(\xi) = 1, \ \ \ \ L_1(\xi) = \xi. &
\end{array}
\end{equation}
The Legendre polynomials form an orthogonal basis for polynomials on $\hat{K}$ since
\begin{equation}
\int_{-1}^{1} L_p(\xi) L_q(\xi) d\xi = \frac{2}{2p+1} \delta_{pq}.
\end{equation}
This orthogonality property of basis functions  implies that the mass matrices are diagonal and gives improved conditioning to the resulting discretization. In particular, we will use   basis functions
\begin{equation}
\phi^m_l(z)= \sqrt{2l+1}L_l\Big(2\frac{z-z_m}{\Delta z_m}\Big),
\end{equation}
which are normalized so that
\begin{equation}
\int_{z_{m-\frac 12}}^{z_{m+\frac 12}} \phi^m_{p}(z) \phi^m_{q}(z) \,dz = \Delta z_m \delta_{pq}.
\end{equation}
Therefore, each component of the solution of  \eqref{lin1d_comp} will be represented as
\begin{equation}
\label{modalexpdg}
q_k(z,t)\approx \sum_{j=0}^{p}  q^{(j)}_{k,m}(t) \phi^m_j(z), \ \ \ \ z\in K_m.
\end{equation}
Integrating \eqref{lin1d_comp} against a test function $\phi^m_i(z)$ on $K_m$
one obtains for $k=1,\dots,d,  $  $i=1,\dots,N $ and
$m=0,\dots,p$

\begin{eqnarray}
\int_{z_{m-\frac 12}}^{z_{m+\frac 12}} \frac{\partial q_k}{\partial t} \phi^m_i(z) \,dz 
&=&
- \sum_{l=1}^d  \int_{z_{m-\frac 12}}^{z_{m+\frac 12}}  a_{kl}   \frac{\partial q_l}{\partial z}   \phi^m_i(z) \,dz
\nonumber\\
 &+& \sum_{l=1}^d \int_{z_{m-\frac 12}}^{z_{m+\frac 12}} b_{kl} q_l(z,t) \phi^m_i(z)  \,dz.
\label{dg_sys_1}
\end{eqnarray}
In the linear, constant coefficient case 
 a full spatial semi-discretization is obtained then by substituting
 \eqref{modalexpdg}, performing standard integration by parts
 and introducing numerical fluxes, see e.g. \cite{cockburn:1989b}. In the variable coefficient or nonlinear
 case, the double integration by part technique proposed in \cite{bassi:1997} (see also \cite{tumolo:2015})
 is employed to handle the non conservative product,
 as well as numerical integration of the resulting integrals by an appropriate quadrature rule.
 
 The discrete equations derived from \eqref{dg_sys_1} will require the knowledge
 of the approximate value of $\lim_{z\rightarrow L^+}q(z,t), $ which will be provided by 
 the semi-infinite Laguerre approximation $\sum_{j=0}^M   q_{l,j}(t).$ Analogously,
 the boundary condition terms $g_l(t) $ in \eqref{modal_sys_2} will be provided by
 the DG approximation
 $\sum_{j=0}^{p}  q^{(j)}_{k,N}(t) \phi^N_j(L).$ In this way, a seamless integration of the two
 approaches is achieved.

As the focus of the present work is on spatial discretization aspects, for the sake of simplicity the numerical tests of the algorithms described above are performed with a simple explicit time discretization. Writing 
the discrete equations  resulting from \eqref{modal_sys_2}, \eqref{dg_sys_1} as a single ODE system

\begin{equation}
 \frac{d}{dt}\mathbf{q}=\mathbf{f}(t, \mathbf{q}(t))\label{eq:q_fq}
\end{equation}
we divide the simulation interval $[0,\,T]$ into a succession of time instants $0=t^0,\,t^1,\ldots,\,t^{N_t}=T$, with $t^{n+1}=t^n+\Delta t$ and $\Delta t$ the time step, assumed constant. Denoting the approximate solution of \eqref{eq:q_fq} with $\mathbf{q}^n\approx\mathbf{q}(t^n)$, we use the third-order Runge-Kutta scheme \cite{hundsdorfer:1995}:
\begin{align}
 \mathbf{q}^{n+1}&=\mathbf{q}^n+\dfrac{\Delta t}{6}\left( \mathbf{K}_1+ \mathbf{K}_2+4 \mathbf{K}_3\right)\\
  \mathbf{K}_1&= \mathbf{f}\left(t^n;\,\mathbf{q}^n\right);\\
  \mathbf{K}_2&= \mathbf{f}\left(t^n+\Delta t;\, \mathbf{q}^n+\Delta t \mathbf{K}_1\right);\\
  \mathbf{K}_3&= \mathbf{f}\left(t^n+\dfrac{\Delta t}{2};\mathbf{q}^n+\dfrac{\Delta t}{4}\left( \mathbf{K}_1+ \mathbf{K}_2\right)\right)
\end{align}
A semi-implicit time discretization stategy will be considered in future multi-dimensional extensions of this paper.
\section{Numerical experiments}
\label{tests}
In this Section we report the results of numerical tests with the coupled DG-Laguerre model described in Section \ref{coupled}. In view of the results of the analysis, we will consider a weak modal discretization based on scaled Laguerre functions with GLR quadrature in the semi-infinite domain. Linear polynomials will be used for the modal DG scheme in the finite domain. 

The coupled scheme will be tested on the case of
the shallow water equations and the suite of tests of \cite{benacchio:2013}. The accuracy of the coupling strategy will be verified by sending signals from either side of the interface between the finite and semi-infinite domain and evaluating the solution by comparison with a modal DG discretization on a single domain of the same extension as the coupled model. As in \cite{benacchio:2013}, a large number
of modes will be used in the semi-infinite part for this first test in order to identify the coupling error.

By using a damping term in the semi-infinite part, the coupled model will be used to simulate an absorbing layer, with the aim to show that outgoing waves can be accurately and efficiently dissipated with minimal reflections into the finite region. To this end the model will initially be run with a long lead time and a non homogeneous wavetrain inflow boundary condition in the finite domain, again comparing the final result with a single-domain discretization. Next, the efficiency of the absorbing layer will be tested by lowering the number of spectral modes and assessing the amplitude of reflected signals in tests with Gaussian initial data. 

%
%
%
The equations governing the small perturbations of a free surface of a nonrotating fluid of height $h$ and velocity $u$ under constant gravity for $z\in\mathbb{R}^+=[0,+\infty]$ are:
\begin{align}
&\dfrac{\partial h}{\partial t}+U\dfrac{\partial h}{\partial z}+H\dfrac{\partial u}{\partial z}+\gamma h=0\label{lsw_h}\\
&\dfrac{\partial u}{\partial t}+g\dfrac{\partial h}{\partial z}+U\dfrac{\partial u}{\partial z}+\gamma u=0\label{lsw_u}
\end{align}
together with initial data and boundary conditions. $H$ and $U$ are respectively the reference height and velocity, and $g$ is the acceleration of gravity. Moreover, a Rayleigh damping term with coefficient $\gamma$ is included. The system \eqref{lsw_h}-\eqref{lsw_u} can be derived from the compressible fluid flow equations with the assumption of small aspect ratio of vertical and horizontal length scales \cite{decoene:2009}.
\subsubsection*{Coupling validation test}
\label{sec:coupl_valid}
\begin{figure}[h]
\centering
\includegraphics[width=.45\textwidth]{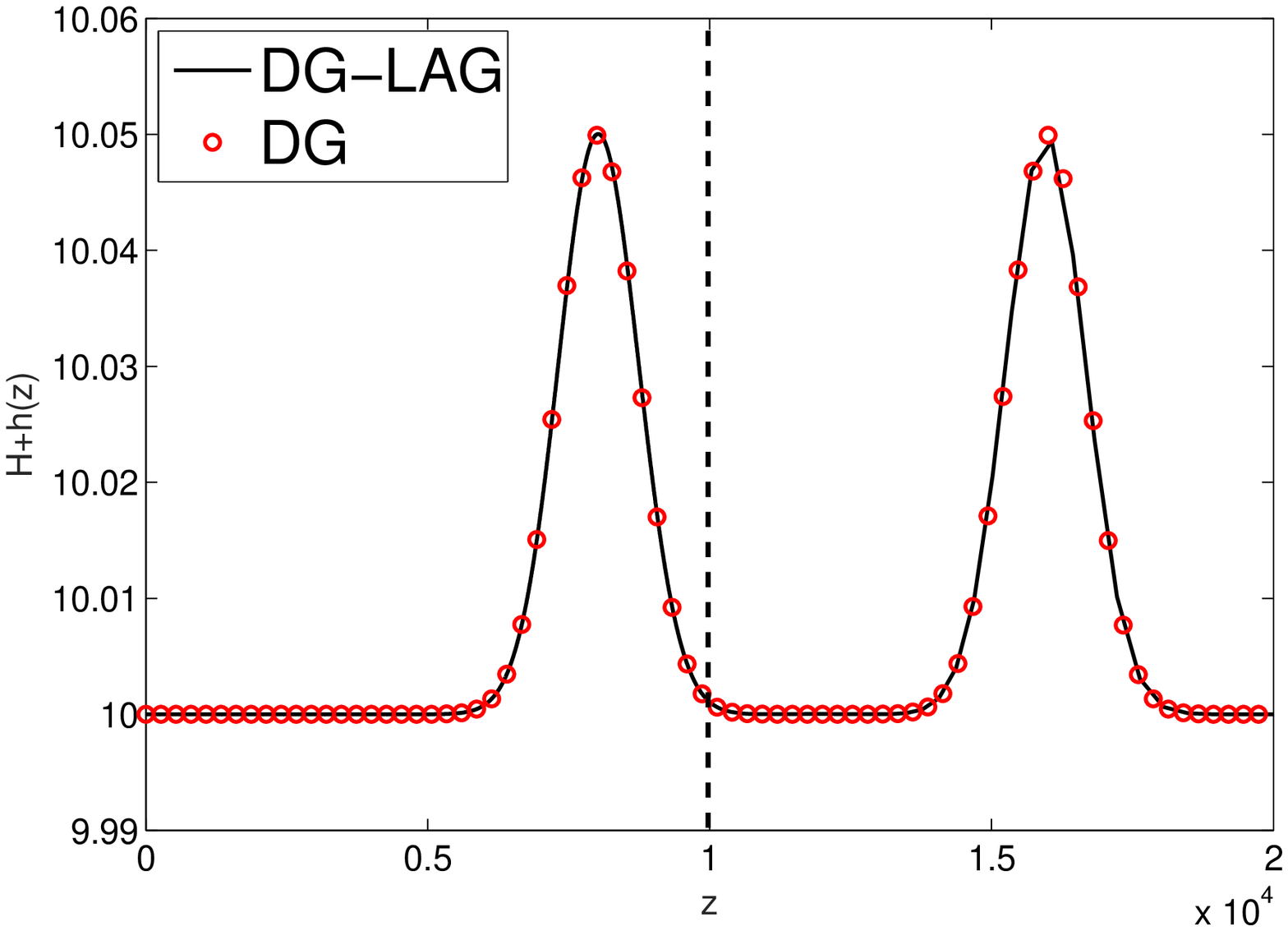}\quad
\includegraphics[width=.45\textwidth]{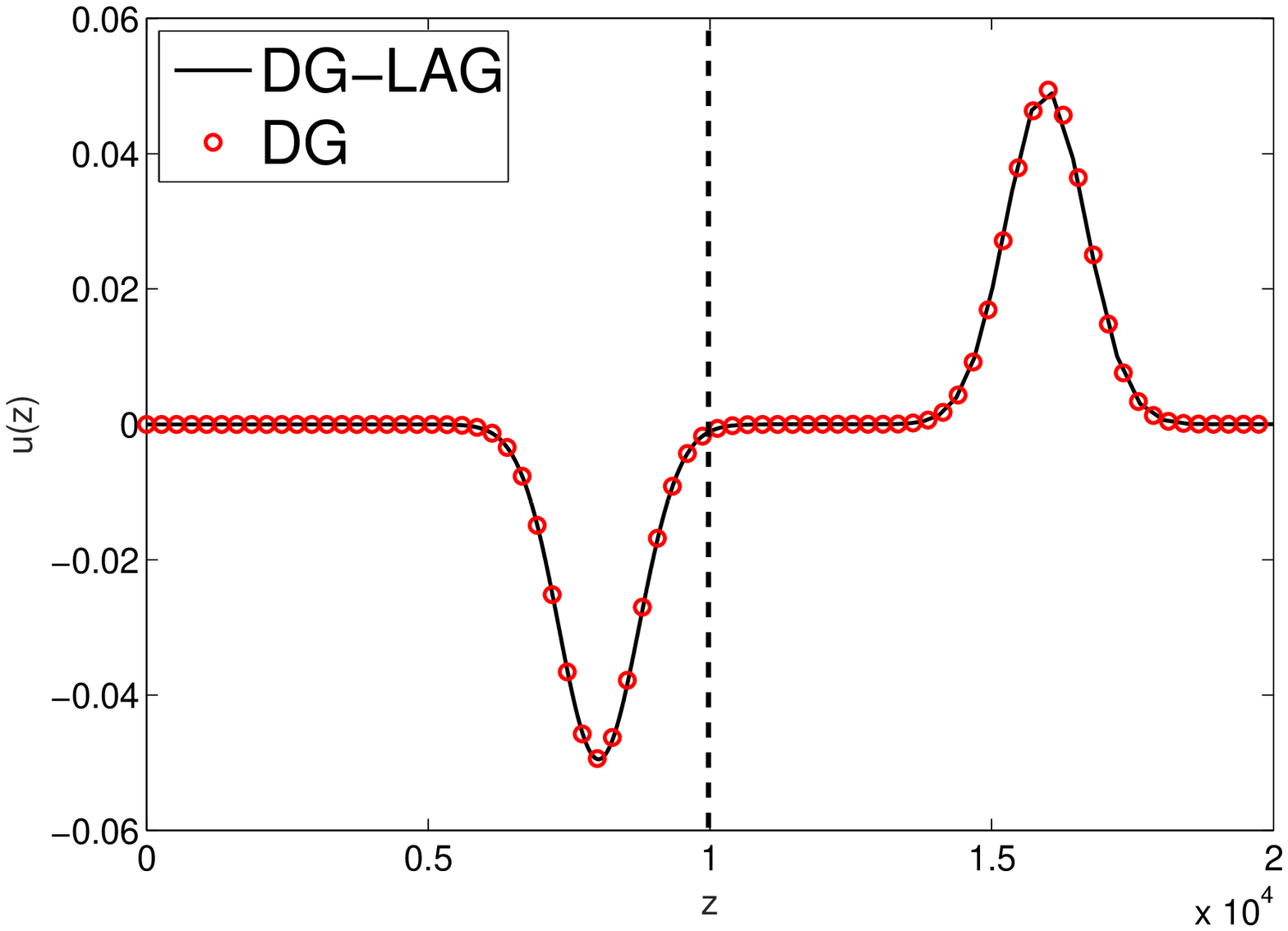}
\caption{Computed final height (left) and velocity (right) by the coupled DG-Laguerre scheme (solid black line)  and single-domain DG scheme (red dots) in the coupling validation, ingoing wave test, $A=0.1\,\textrm{m},\,\sigma=1000\,\textrm{m}$.}
\label{fig:coupl_valid_A0.1_sigma1000}
\end{figure}

The first test uses a motionless Gaussian initial distribution for the height, $h(x, 0)=h_1\exp{[(x-x_0)/\sigma]}$, where we set $h_1=0.1\,\textrm{m}$ and consider two choices of $\sigma=500,\,1000\,\textrm{m}$. The number of modes in the semi-infinite domain is $N=180$, and $\beta=1/400$, so that the last GLR point is at $x_N=2.86\textrm{E}05\,\textrm{m}$. The absorbing layer is not active, $\gamma=0$. We consider a $10000\,\textrm{m}$-wide finite domain, with $N_x=1250$ cells. We then distinguish the cases of an ingoing wave, $x_0=12000\,\textrm{m}$, and an outgoing wave, $x_0=5000\,\textrm{m}$. 
%
%
The initial datum splits in two crests and the leftmost one travels through the interface between the finite and semi-infinite domain. The solution of the coupled model is compared on the finite domain with a reference solution obtained with a full DG discretization over $20\,\textrm{km}$ using the same spatial resolution and number of time steps $N_t=2200$ (Figure \ref{fig:coupl_valid_A0.1_sigma1000}). As a large number of modes is used in the semi-infinite domain, the relative errors in the finite domain $[0,\,10]\,\textrm{km}$ between the multidomain and the single domain schemes (Table \ref{tab:err_coupl_valid}) identify the residual perturbations coming from the coupling approximation at the interface. Similarly to \cite{benacchio:2013}, the errors are at most around a few percent. 
\begin{table}[h]
\centering
\scriptsize
\caption{Computed 1-norm ($\mathcal{E}_1$), 2-norm ($\mathcal{E}_2$), and maximum ($\mathcal{E}_\infty$) errors for the coupling validation experiment, ingoing wave (top four rows, relative errors) and
outgoing wave (bottom four rows, absolute errors) cases.}
\begin{tabularx}{\textwidth}{cccXXXXXX}\toprule\midrule
$x_0$  & $\tilde{h}$ & $\sigma$ & $\mathcal{E}_1(\eta)$ & $\mathcal{E}_1(\mathbf{u})$ & $\mathcal{E}_2(\eta)$ & $\mathcal{E}_2(\mathbf{u})$ &   $\mathcal{E}_\infty(\eta)$ & $\mathcal{E}_\infty(\mathbf{u})$ \\\midrule
12000 & 0.1 & 1000  & 7.37E-03 & 7.37E-03 & 8.49E-03 & 8.48E-03 & 1.10E-02 & 1.10E-02\\
12000 & 0.1 & 500  & 1.58E-02 & 1.57E-02 & 1.70E-02 & 1.70E-02 & 1.96E-02 & 1.96E-02\\
12000 & 0.5 & 1000  & 3.67E-02 & 3.65E-02 & 4.14E-02 & 4.12E-02 & 5.11E-02 & 5.09E-02\\
12000 & 0.5 & 500  & 7.81E-02 & 7.78E-02 & 8.33E-02 & 8.30E-02 & 8.57E-02 & 8.54E-02\\\midrule
5000 & 0.1 & 1000  & 3.88E-06 & 3.84E-06 & 9.23E-06 & 9.14E-06 & 3.11E-05 & 3.08E-05\\
5000 & 0.1 & 500  & 1.94E-06 & 1.92E-06 & 6.52E-06 & 6.46E-06 & 3.11E-05 & 3.08E-05\\
5000 & 0.5 & 1000  & 9.36E-05 & 9.27E-05 & 2.23E-04 & 2.21E-04 & 7.58E-04 & 7.50E-04\\
5000 & 0.5 & 500  & 4.68E-05 & 4.64E-05 & 1.58E-04 & 1.56E-04 & 7.58E-04 & 7.50E-04\\
\bottomrule\end{tabularx}
\label{tab:err_coupl_valid}
\end{table}
%
%
In a second test, the initial datum is placed inside the finite domain and run for
$N_t=8400$ time steps until a final time $T=1000\,\textrm{s}$ when all perturbations have left the finite domain. Since the solution in the finite domain is the absence of perturbations in the finite domain, absolute errors  with respect to the reference are computed in this case (bottom rows in Table \ref{tab:err_coupl_valid}). The residual perturbations in the finite domain have negligible amplitude. While a one-to-one comparison with 
\cite{benacchio:2013}, which coupled two fully nonlinear discretizations, is not possible here, the obtained errors appear small enough to test the coupled DG-Laguerre scheme to simulate absorbing layers. Notice that here a different (and more accurate) time discretization method is employed with respect to \cite{benacchio:2013}.
%
%
Next, we test the coupled scheme with $\gamma\neq0$ in the semi-infinite domain. As in \cite{benacchio:2013}, the damping coefficient will have the sigmoid-like functional form:
\begin{equation}
 \gamma(x)=\dfrac{\Delta\gamma}{1+\exp\left(\dfrac{\alpha L_0-x}\sigma\right)}
\end{equation}
where $\Delta\gamma$ is the sigmoid amplitude, $L_0$ a length scale, typically the absorbing layer thickness, $\alpha$ the position of the sigmoid inside the absorbing layer, and $\sigma$ the sigmoid steepness,   see also Figure 8 in \cite{benacchio:2013}, to which we refer for a complete definition of the parameter values employed.

\begin{figure}[h]
\centering
\includegraphics[width=.45\textwidth]{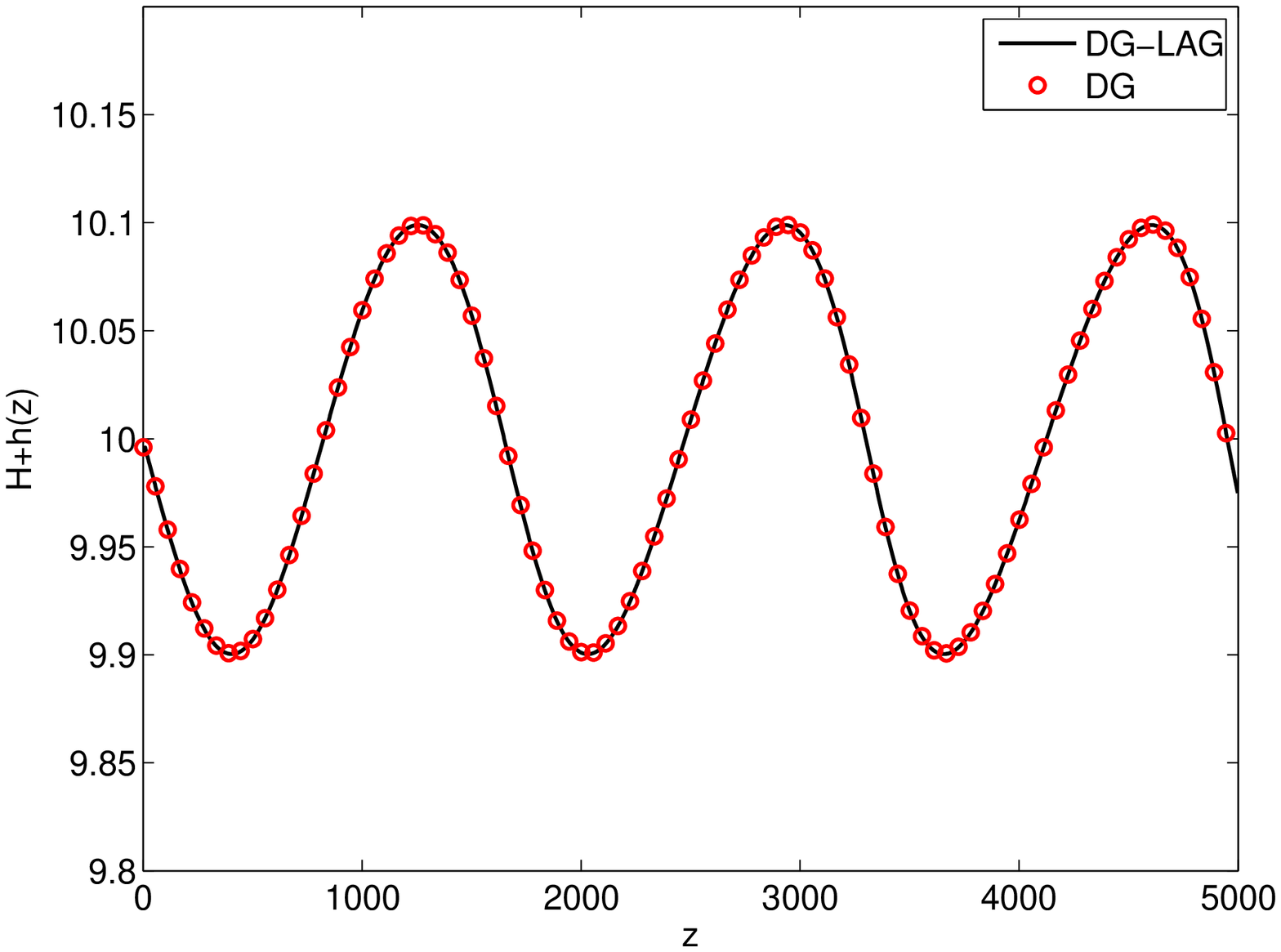}
\includegraphics[width=.45\textwidth]{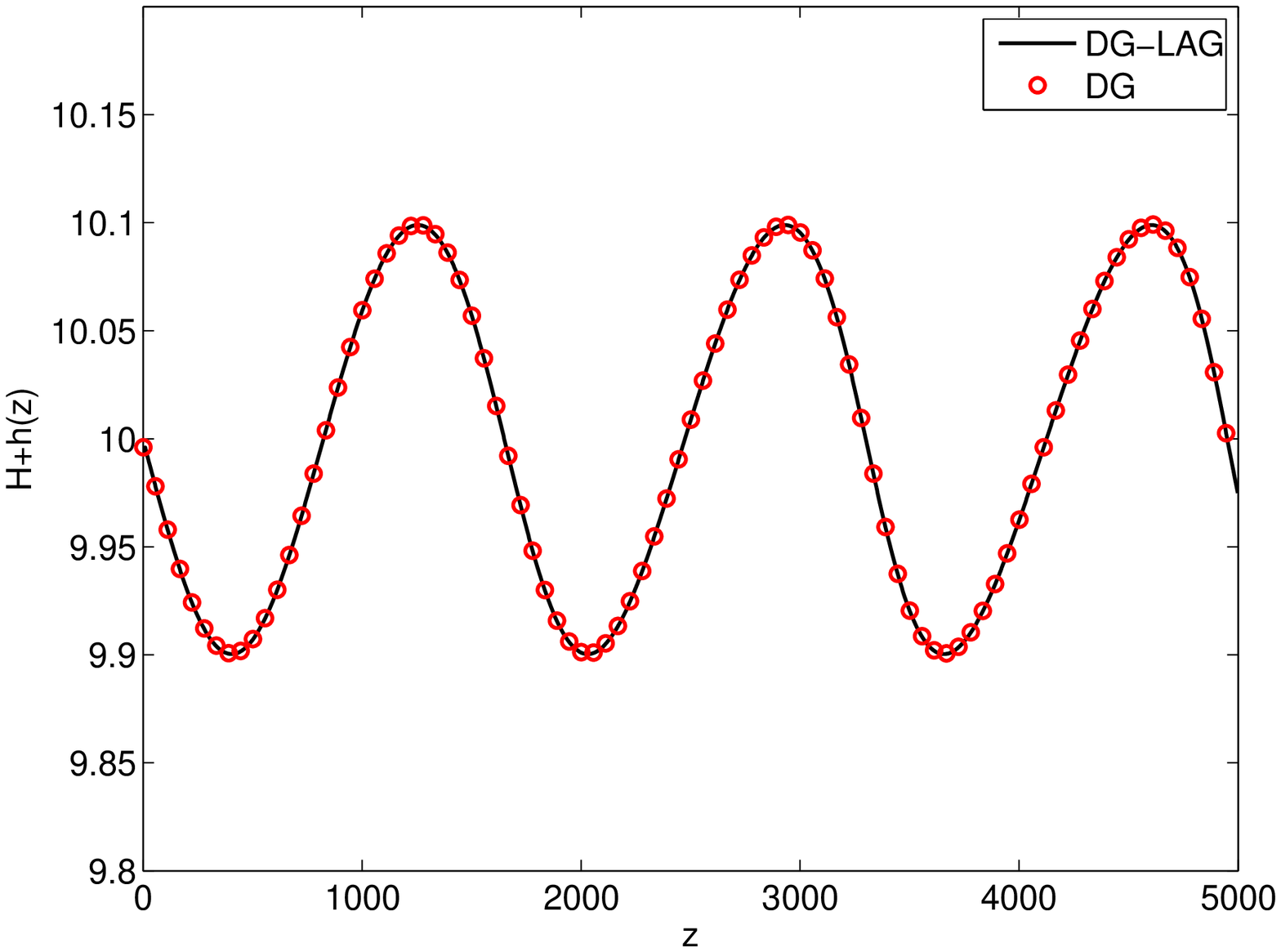}
\caption{Computed total depth after $T=5000\,\textrm{s}$ in the wavetrain test for the coupled DG-Laguerre (solid line) and single-domain DG model (red circles), with 30 (left) and 15 (right) semi-infinite nodes, $A=0.05\,\textrm{m},\, k=30\,\textrm{m$^{-1}$}$.}
\label{fig:wavetrainH_k30}
\end{figure}

\begin{figure}[h]
\centering
\includegraphics[width=.45\textwidth]{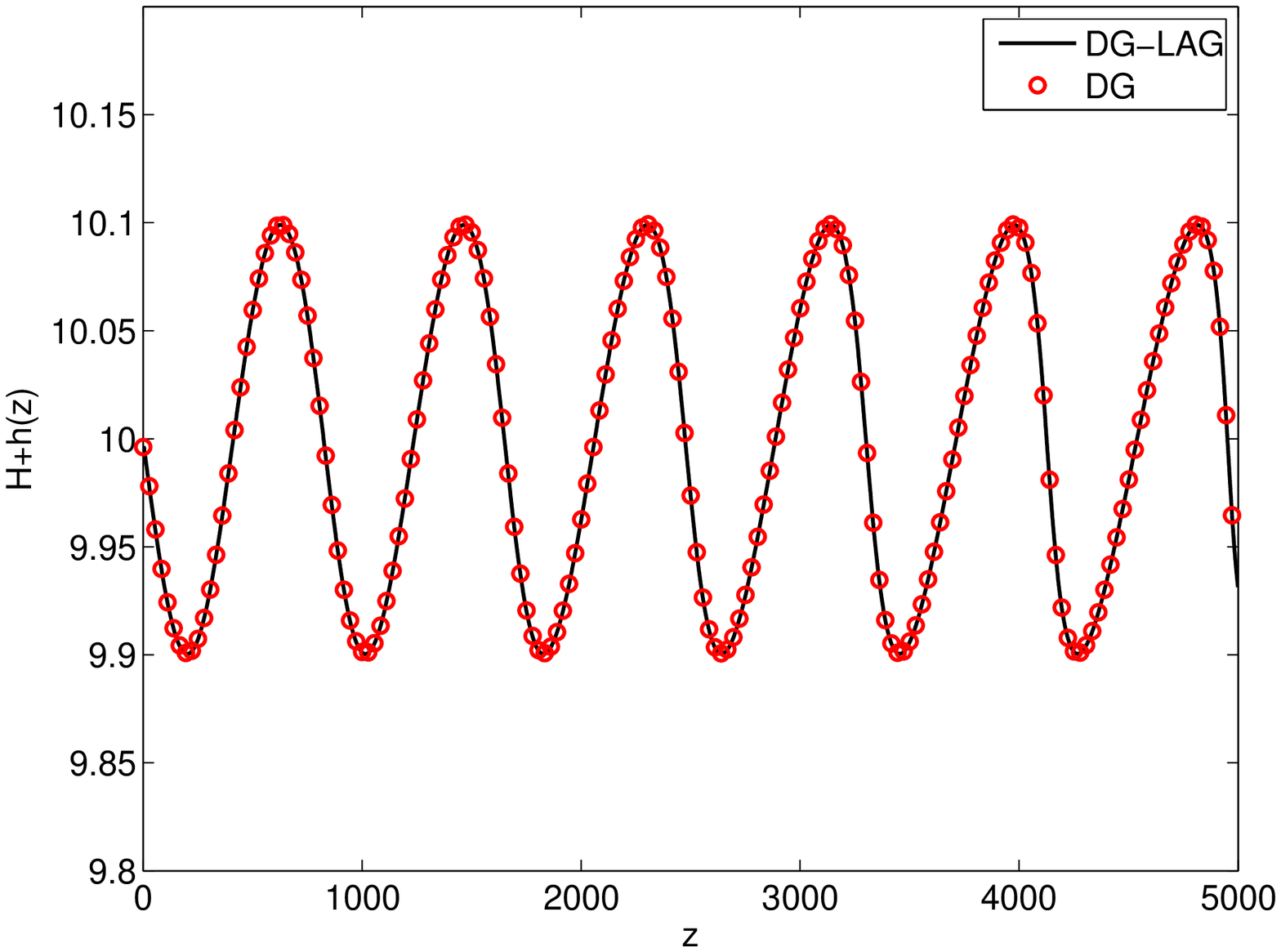}\quad
\includegraphics[width=.45\textwidth]{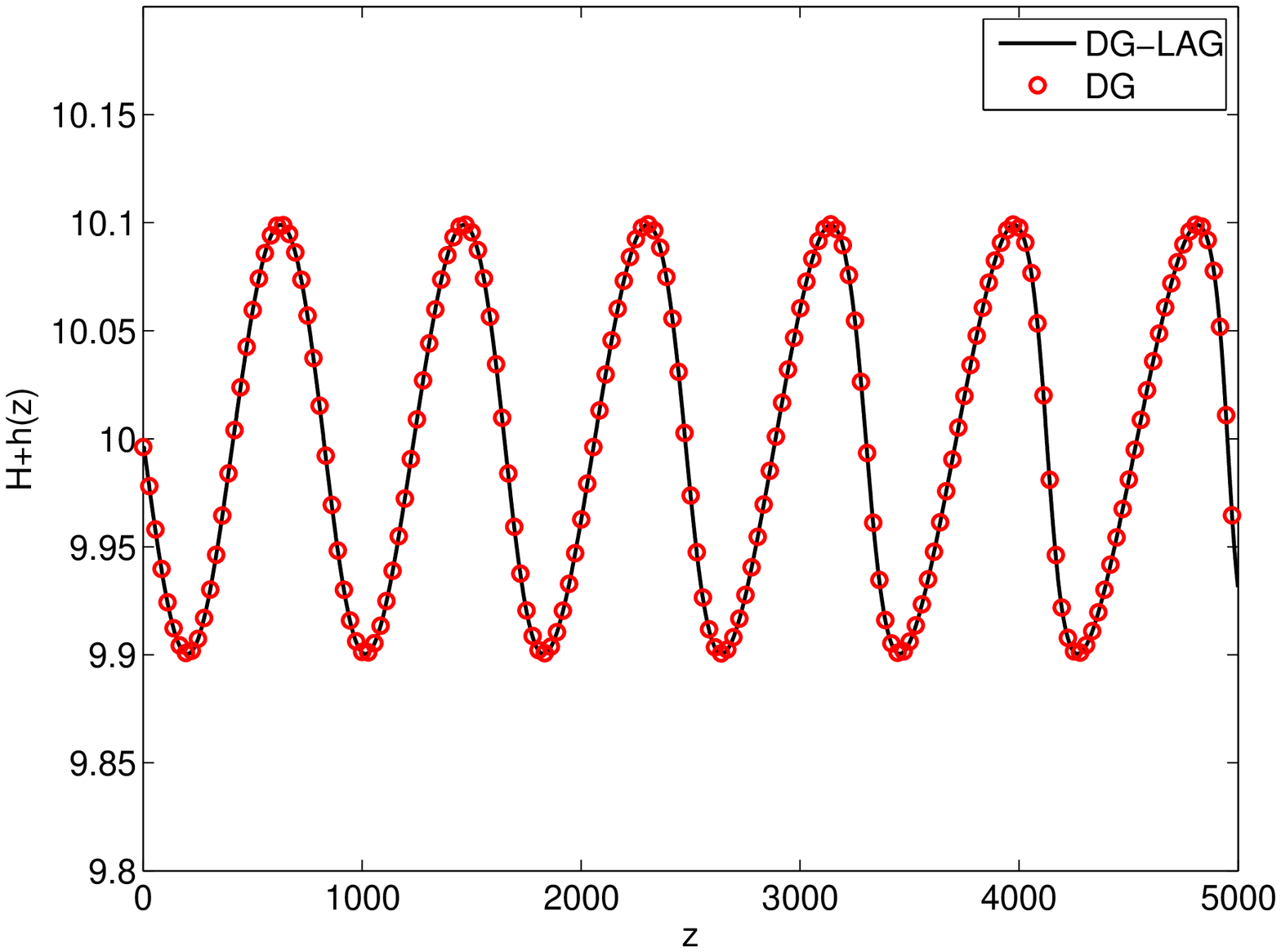}
\caption{Computed total depth after $T=5000\,\textrm{s}$ in the wavetrain test for the coupled DG-Laguerre (solid line) and single-domain DG model (red circles), with 30 (left) and 15 (right) semi-infinite nodes, $A=0.05\,\textrm{m},\, k=60\,\textrm{m$^{-1}$}$.}
\label{fig:wavetrainH_k60}
\end{figure}

A Dirichlet boundary condition on the velocity is imposed in $x=0$: 
$u(0,\,t)=A\sin({2\pi k/T})$. The initial condition is motionless, with $h(0)=0$. With these choices, a train of waves with uniform amplitude is generated and propagated through the $5\,\textrm{km}$-wide finite domain. Once the waves cross the interface, they are damped in the semi-infinite domain that acts as an absorbing layer. In this context, excessive spurious signals coming from the interface would propagate and pollute the wavetrain in the finite domain. We compare again the solution of the coupled model with the one obtained with a DG scheme over a larger domain, with $N=30$ or $N=15$ modes in the semi-infinite domain and a longer final time $T=5000\,\textrm{s}$, with two choices for the wavenumber, $k=30,\,15\,\textrm{m$^{-1}$}$ and three choices for the amplitude, $A=0.025,\,0.05,\,0.1\,\textrm{m}$. The case provides a simplified representation of gravity wave propagation in the atmosphere. The obtained relative errors on height and velocity as well as the energy errors:
\begin{equation}
 \mathcal{E}_{EN}=\dfrac1N_x\sum_i\dfrac12\left[g(h_i-[h_{\textrm{ref}}]_i)^2+H(u_i-[u_{\textrm{ref}}]_i)^2\right].
\end{equation}
are low for all configurations (Figures \ref{fig:wavetrainH_k30} and \ref{fig:wavetrainH_k60} and Tables \ref{tab:wave_train_30nodes} and \ref{tab:wave_train_15nodes}) and in line with the results in \cite{benacchio:2013}. 
For the evaluation of these results, we note that the Dirichlet condition was imposed on the discharge in \cite{benacchio:2013} so results obtained with, for instance, $A=0.2\,\textrm{m}$ in the present paper should be compared to results with $A(H+\eta)\approx2\,\textrm{m}$ of \cite{benacchio:2013}.

\begin{table}[h]
\centering
\scriptsize
\caption{Relative mean square root and maximum errors for elevation $\eta$ and velocity $\mathbf{u}$ for the wavetrain test, 30 nodes in the semi-infinite domain. See text for other parameters.}
\label{tab:wave_train_30nodes}
\begin{tabularx}{\textwidth}{ccccXXXXX}\toprule\midrule
$A$   & $k$ & $N_x$ & $\beta$ & $\mathcal{E}_2^{\textrm{rel}}(\eta)$ & $\mathcal{E}_\infty^{\textrm{rel}}(\eta)$ & $\mathcal{E}_2^{\textrm{rel}}(\mathbf{u})$ & $\mathcal{E}_\infty^{\textrm{rel}}(\mathbf{u})$ & $\mathcal{E}_{\textrm{EN}}$\\\midrule
0.025 & 30  & 600    & 0.0143 & 3.84E-06 & 4.40E-04 & 6.23E-06 & 6.10E-04 & 8.38E-09\\
0.025 & 60  & 1200    & 0.0286 & 3.89E-06 & 6.50E-04 & 6.03E-06 & 8.10E-04 & 9.94E-09\\
0.05 & 30  & 600    & 0.0143 & 1.67E-05 & 8.88E-04 & 2.72E-05 & 1.46E-03 & 1.55E-07\\
0.05 & 60  & 1200    & 0.0286 & 2.54E-05 & 1.19E-03 & 5.41E-05 & 2.14E-03 & 3.51E-07\\
0.1 & 30  & 600    & 0.0143 & 1.01E-04 & 1.98E-03 & 1.66E-04 & 5.22E-03 & 5.41E-06\\
0.1 & 60  & 1200    & 0.0286 & 1.32E-04 & 5.08E-03 & 1.50E-03 & 6.71E-02 & 1.11E-05\\
\bottomrule\end{tabularx}
\end{table}
\begin{table}
\centering
\scriptsize
\caption{Relative mean square root and maximum errors for elevation $\eta$ and velocity $\mathbf{u}$ for the wavetrain test, 15 nodes in the semi-infinite domain. See text for other parameters.}
\begin{tabularx}{\textwidth}{ccccXXXXX}\toprule\midrule
$A$   & $k$ & $N_x$ & $\beta$ & $\mathcal{E}_2^{\textrm{rel}}(\eta)$ & $\mathcal{E}_\infty^{\textrm{rel}}(\eta)$ & $\mathcal{E}_2^{\textrm{rel}}(\mathbf{u})$ & $\mathcal{E}_\infty^{\textrm{rel}}(\mathbf{u})$ & $\mathcal{E}_{\textrm{EN}}$\\\midrule
0.025 & 30  & 600    & 0.0286 & 3.80E-06 & 4.37E-04 & 6.17E-06 & 6.07E-04 & 8.22E-09\\
0.025 & 60  & 1200    & 0.0571 & 6.84E-06 & 5.34E-04 & 2.27E-05 & 2.64E-03 & 2.47E-08\\
0.05 & 30  & 600    & 0.0286 & 1.65E-05 & 8.83E-04 & 2.70E-05 & 1.45E-03 & 1.53E-07\\
0.05 & 60  & 1200    & 0.0571 & 3.68E-05 & 1.38E-03 & 1.13E-04 & 7.21E-03 & 7.09E-07\\
0.1 & 30  & 600    & 0.0286 & 1.01E-04 & 1.97E-03 & 1.66E-04 & 5.21E-03 & 5.36E-06\\
0.1 & 60  & 1200    & 0.0571 & 1.98E-04 & 9.59E-03 & 2.49E-03 & 1.12E-01 & 2.82E-05\\
\bottomrule\end{tabularx}
\label{tab:wave_train_15nodes}
\end{table}
%

%
%
In the final test of the suite of \cite{benacchio:2013}, we evaluate the performance of the coupled model in absorbing signals leaving the finite domain. By using a small number of nodes in the semi-infinite domain, the objective is to show that outgoing waves can be efficiently damped and that reflected signals at the interface have low amplitude. As in Section \ref{sec:coupl_valid} above, we consider a $D=10\,\textrm{km}$-wide domain and a Gaussian perturbation centred in $x_0=7500\,\textrm{m}$ with amplitude $\sigma=500\,\textrm{m}$. In this case we run the coupled model until the final time $T=D/(2\sqrt{gH})$, at which reflected perturbations will have come back to $x_0$. We consider the same performance measures as in \cite{benacchio:2013}: residual maximum values for free-surface elevation and velocity, mean square root energy error, and reflection ratio, defined as:
\begin{equation}
\rho=\sqrt{\dfrac{\mathcal{E}_{EN}(T)}{\mathcal{E}^W_{EN}(T)}}
\end{equation}
 where $\mathcal{E}^W_{EN}$ denotes the energy error obtained with a solid wall boundary condition at the interface. Results with the DG-Laguerre coupled model are in line or lower than the ones obtained in \cite{benacchio:2013} with a finite volume discretization in the finite domain (Tables \ref{tab:abs_layer_eff1} and \ref{tab:abs_layer_eff2}).  
\begin{table}[h]
\centering
\scriptsize
\caption{Maximum residual elevation, velocity, mean square root energy error and reflection ratios for the single Gaussian perturbation test, coupled DG/Laguerre scheme.}
\begin{tabularx}{.7\textwidth}{ccccccc}\toprule\midrule
$N$  & $N_x$ & $T/\Delta t$ & $\|\eta(T)\|_\infty$ & $\|u(T)\|_\infty$ & $\mathcal{E}_{\textrm{EN}}$ & $\rho$ \\\midrule
40 & 400 & 600 & 2.92E-03 & 2.90E-03 & 7.23E-06 & 4.57E-03 \\
30 & 300 & 450 & 2.91E-03 & 2.88E-03 & 7.17E-06 & 4.55E-03 \\
20 & 200 & 300 & 3.00E-03 & 2.97E-03 & 7.27E-06 & 4.58E-03 \\
10 & 100 & 150 & 3.13E-03 & 3.11E-03 & 8.16E-06 & 4.88E-03 \\
\bottomrule\end{tabularx}
\label{tab:abs_layer_eff1}
\end{table}
\begin{table}[h]
\centering
\scriptsize
\caption{Further reflection ratios for the single Gaussian perturbation test, coupled DG/Laguerre scheme.}
\begin{tabularx}{.4\textwidth}{ccccc}\toprule\midrule
$N$  & $N_x$ & $T/\Delta t$ & $\beta$ & $\rho$ \\\midrule
40 & 400 & 600 & 1/280 & 4.57E-03 \\
30 & 400 & 600 & 1/210 & 4.57E-03 \\
20 & 400 & 600 & 1/145 & 4.56E-03 \\
10 & 400 & 600 & 1/75 & 4.56E-03 \\
5 & 400 & 600 & 1/40 & 4.05E-03 \\\midrule
30 & 300 & 450 & 1/280 & 4.55E-03 \\
20 & 300 & 450 & 1/190 & 4.54E-03 \\
10 & 300 & 450 & 1/110 & 4.52E-03 \\
5 & 300 & 450 & 1/40 & 4.12E-03 \\\midrule
20 & 200 & 300 & 1/280 & 4.58E-03 \\
10 & 200 & 300 & 1/120 & 4.52E-03 \\
5 & 200 & 300 & 1/55 & 4.21E-03 \\\midrule
13 & 110 & 150 & 1/270 & 4.80E-03 \\
5 & 110 & 150 & 1/250 & 5.03E-03 \\
\bottomrule\end{tabularx}
\label{tab:abs_layer_eff2}
\end{table}

\section{Conclusions and perspectives}
\label{conclu}
We analyzed the stability of Laguerre spectral discretizations of hyperbolic problems on semi-infinite domains, by computing the spectra for the strong and weak, nodal and modal discretization of the 
linear advection equation based on Laguerre polynomials and functions with a range of options for boundary conditions and quadrature rules. Discretizations using 
Gauss-Laguerre-Radau quadrature nodes and scaled Laguerre functions
 were found to give the best results in terms of stability, while the analysis rules out the use of Laguerre polynomials due to their poor stability properties. To the best of our knowledge, this is the first analysis of this kind in the literature on spectral methods for fluid dynamics.

The theory was then extended to systems of equations, and a modal semi-infinite discretization was coupled with a discontinuous finite element scheme in a finite domain. The resulting scheme was tested on the
propagation of single wave and wavetrains, with outcomes on validation and absorbing layer efficiency tests in line with the findings in \cite{benacchio:2013}. Therefore, the scaled Laguerre approach
in the unbounded region was effectively shown to work independently of the form of the discretization in the finite domain. 


The results suggest the potential of the Laguerre spectral discretization as a flexible and independent tool to accurately simulate wave propagation on unbounded domains, as well as to reduce the computational
cost of absorbing unwanted perturbations without spurious reflections. In particular, extensions to nonlinear systems can be considered in multiple space dimensions via tensor-product
approaches in the unbounded domain and continuous, discontinuous, or mixed finite element discretizations
in the bounded domain \cite{cotter:2012}. Due to the variety of interacting solutions in full-fledged numerical weather prediction systems, the flexibility of the scaled spectral absorbing layer could
be productively employed to lower the cost of standard approaches currently deployed in operational models. 


%
\section*{Acknowledgements} 
We kindly acknowledge useful discussions on the topics of this paper with Vadym Aizinger.
The contribution of T.B. is covered by {\textcopyright Crown Copyright 2018}. L.B. gratefully acknowledges
the financial support from the Met Office for a visit in May 2015, in which the work on this paper was started.
\bibliographystyle{plain}
\bibliography{laguerre}

\end{document}